\documentclass[10pt]{article}
\usepackage{amssymb}

\usepackage{amsthm}

\begin{document}

\title{SOME REMARKS ON MIXED PROBLEMS}
\author{Tove Dahn}
\maketitle

\section{Introduction}
The model we start with can be compared with a lifting principle, $F(\gamma)(\zeta)=f(\zeta)$, $\zeta \in \mathbf{R}^{n}$
where $\widehat{g}(\zeta)=f(\zeta)$ and $\widehat{E}=F$. When $\gamma$ polynomial in an ideal $(I)$, 
we assume thus that we have existence of $F$ analytic, such that $E * \gamma(D) \delta=g$. 
More generally we can assume existence of $E \in \mathcal{D}_{L^{1}}'$,such that $E(\gamma)=g$. 
We will restrict our attention to symbols $f$ with $\log \mid f \mid \in L^{1}$ and for a very regular boundary
(\cite{Dahn13}), we argue that the model can be represented using polynomial $\gamma$.
By further using the moment problem, we can continue $F$ to $C$, 
given that the measure to $(I)^{\bot}$ is of bounded variation. If further for instance the measure is reduced, 
we can continue to $\mathcal{E}'$. 

We consider it as necessary for a global model to have a representation of the complementary ideal.
Assume the analytic symbols have a decomposition $(I) \oplus (I)^{\bot}$, we use the moment problem 
to generalize this orthogonality (\cite{Riesz56}). 
Assume that we have a continuous mapping $(I) \ni \gamma_{1} \rightarrow \gamma_{2} \in (I^{\bot})$. 

A global model is invariant for change of local coordinates. We discuss instead a model invariant 
for change of pairs of coordinates $H(\gamma)=H_{1}H_{2}(\gamma)$. Note that if $H$ polynomial, 
we do not necessarily have $H_{1},H_{2}$ polynomials. We discuss $F(\gamma_{1} \rightarrow \gamma_{2})(\zeta)$, 
where $\gamma_{1}$ is hyperbolic and $\gamma_{2}$ partially hypoelliptic. 

We assume an invariance principle for movements in $(I)$, so that any movement of $(I)$ has an exact
correspondent movement in $(I)^{\bot}$ and conversely. When the movement $V$ has an analytic representation, we have a continuous mapping $V \rightarrow \tilde{V}$ into
the domain for $\zeta$. The main result that we will discuss in this paper is:
\newtheorem*{res15}{Main result}
\begin{res15}
 Assume $\gamma_{1} \in (I_{1})$ an analytic symbol and consider a continuation $(I_{1}) \ni \gamma_{1} \rightarrow \gamma_{\delta} \in (I_{1})^{\bot}$,
 for a parameter $\delta$. We assume there is a set $\Delta$ (lineality), where $\gamma_{1}=\gamma_{\delta}$.  
 We then have, given $f(\zeta)$ analytic, a solution $F_{\delta}$, such that $F_{\delta}(\gamma_{\delta})(\zeta)=f(\zeta)$. 
 We assume the continuation such that 
 $\frac{d F_{\delta}}{d x_{1}}/\frac{d F_{\delta}}{d y_{1}}=\frac{d F_{\delta}}{d x_{\delta}}/\frac{d F_{\delta}}{d y_{\delta}}$, under
 the following conditions on the boundary: The boundary is given by the first surfaces $S$ to $f(\zeta)$,
 that can be reached in at least one of three ways.
 \begin{itemize}
  \item[1] $S$ can be reached by $\Delta$
  \item[2] $S$ can be reached through a movement $V^{\bot} \leftarrow U^{\bot}$, such that
  $\tilde{V^{\bot}}$ and $\tilde{V}$ do not have points in common outside $S$.
  \item[3] $S$ can be reached through an algebraic trajectory.
 \end{itemize}
  \end{res15}
Note that the mixed problem usually is given as $F \gamma_{\delta}=f_{\delta}$. But if we assume $\gamma_{1}$
has a fundamental solution $E_{1}$ and ${}^{t} E_{1} {}^{t }\gamma_{\delta} F=F \gamma_{\delta} E_{1}$,
then ${}^{t} E_{1} {}^{t} \gamma_{\delta} F \gamma_{1} = f_{\delta}$, where $f_{\delta}={}^{t} E_{1} {}^{t} \gamma_{\delta} f$

\section{The invariance principle}

We start with (\cite{Riesz43}) a movement on the hyperboloid $L(Ux,Uy)=L(x,y)$, that is ${}^{t}U \sim U^{-1}$ 
with respect to a lorentz metric form $L$.
The reflection invariance is defined by $L(Ux,y)=L(x,y)$ implies $y=0$, in this case $U \sim I$. 
Particularly, we define the light cone $L(x,y)=0$, that is $x \bot y$ with respect to $L$.

Parallel to this we discuss a en movement with respect to euclidean metric $E(V x,V y)=E(x,y)$ that is ${}^{t}V \sim V^{-1}$
with respect to $E$ and we have an axes for invariance $E(V x,y)=E(x,y)$ implies $y=0$ and we conclude $V \sim I$. 
When the movements are reflections, they are involutive.

The idea behind the model 
is that an hyperbolic movement exactly corresponds to an euclidean translation. In the  non- euclidean plane 
we do not have any proper translation, but given two points $p,q$ there is a unique hyperbolic rotation
that maps $p$ onto $q$ and for which ``the line'' through the points is a trajectory. In this manner
the rotation axes is
conjugated to the plane through the origo and the points $p,q$. (\cite{Riesz43}).

Consider now $\gamma_{1} \in (I_{1})$ where orthogonality is relative $L$ and simultaneously $\gamma_{2} \in (I_{1})^{\bot}$
with orthogonality relative $E$. Using Radon Nikodym's theorem, we will continuously continue $\gamma_{1} \rightarrow \gamma_{2}$.
Assume $< \gamma_{1},\gamma_{2}>=0$ and $U \gamma_{1} \rightarrow I \gamma_{1}$ through a closed sequence
(continuous), say $M$. We can argue that $M={}^{\circ} ( M^{\circ} )$, that is if the condition
$< (U-I) \gamma_{1},(V-I) \gamma_{2}>=0$ implies $V \rightarrow I$ and conversely, when $V=I$, we have that $U=I$ and we have a geometric invariance principle.
Assume $J : \gamma_{1} \rightarrow \gamma_{2}$ continuous and $VJ=JU$, then
using the relation above for $M$ we have motivated an invariance principle.
We are here assuming
$< U \gamma_{1}, V \gamma_{2}>=0$ and $< (U-I) \gamma_{1},\gamma_{2}>=< \gamma_{1},(V-I) \gamma_{2}>=0$ on $\Delta$.

Now for the particular movements, we have
$\tau J = J h$, $e J= J \tilde{e}$ where $\tau$ is translation, $h$ a hyperbolic movement and $e,\tilde{e}$
rotation. Characteristic for
parabolic movements is a constant (euclidean) distance to the light cone (rotation axes).
Characteristic for elliptic and hyperbolic movements is that the rotation axes $R$ are one sided, that is
$\mid x \mid < \mid y \mid$ or $\mid x \mid> \mid y \mid$. Thus if we consider $\eta(x)=y/x$ we have 
that the axes $R_{1} \rightarrow \mid \eta \mid < 1$ (hyperbolic) and $R_{2} \rightarrow \mid \eta \mid > 1$ (elliptic)
and $R_{3} \rightarrow \mid \eta \mid=1$ (parabolic). We can as usual map $\mid \eta \mid < 1$
on a half plane, in this way we can consider movement as ``one sided''. In particular the  hyperbolic and elliptic movements map
half planes on to half planes. 
Elliptic rotation can be represented through $e^{i<a \cdot x,a^{*} \cdot x^{*}>}$, where $a$ is a 
scalar vector and where $a^{*}$, is defined such that $< a \cdot x, a^{*} x^{*}>=< x,x^{*} >$, 
that is elliptic rotation is immediately represented in euclidean geometry. Using the Fourier Borel 
transform, we can further associate elliptic rotation to translation of the symbol. 

Assume the invariance principle $JU=VJ$ and $x \rightarrow x^{*}$ according to Legendre (reciprocal polars)
and note that $< x,x^{*}>=< x^{*},x>$ implies a normal transform. 
We look for a continuous mapping $V^{\bot} \rightarrow V^{*}$. Assume $V_{e}$ corresponds to 
rotational movement and let $\tilde{V_{e}}$ denote the inverse mapping, that is $V_{e}f \rightarrow \tilde{V_{e}} \zeta$, 
we then have  $V_{e}^{\bot} \ni x \rightarrow x^{*} \in V_{h}$ using Legendre,
that is $(V_{h})^{*} \rightarrow (V_{e}^{\bot})$ (reflexivity). In the same manner  
$(V_{e})^{*} \simeq (V_{h})^{\bot}$ and $(V_{p})^{*} \simeq (V_{p}^{\bot}) \simeq (V_{p})$. 
Formally $(e^{i <x, \cdot > + v(x)} f)^{\vee}= \tau \widehat{f}(x^{*})$. Using Parseval we have equivalent 
sets in $L^{2}$. If we define a regular approximation through $H(U \varphi) \equiv 0$, if we assume 
$\frac{d}{d x}H \neq 0$ and $H$ analytic, we have continuous induced relations for the inverse movements $\tilde{V}$ 

For simplicity we will in this article consider
one movement at a time, when we apply Radon-Nikodym's theorem and since the movements are not 
dependent on sign of $L$, is sufficient to consider positive linear functionals.

Assume $\Delta$ a domain for the movement $I$ and $\gamma_{1}=\gamma_{2}$ on $\Delta$. Assume $(U - I) \gamma_{1}=0$ with respect to $L$
and $(V-I) \gamma_{2}=0$ with respect to $E$ and $F(J U \gamma_{1})=F(e^{v} V \gamma_{2})$ with $v \in L^{1}$
where $v = 0$ on $\Delta$. Note that the last condition is dependent on the movement.
Assume $w(\eta)=w(\frac{y}{x})=v(x,y)$
and assume $x \frac{d \eta}{d x}=-\frac{d y}{d x} - \frac{x}{y}$. Further let $d u=\frac{1}{x} d x - \frac{1}{y} d y$
then we have that $-d \eta = d u^{*}$, the harmonic conjugate. 
As $w=0$, over $\Delta $ we have for a trajectory to $\eta$, that analyticity is preserved under the condition on vanishing flux $\int_{S} d u^{*}=0$
(a transmission property). 

Note (\cite{Bendixsson01})
that for $X=d x/ d t,Y=d y/d t$, $(x X + y Y)/(x Y - y X)$ is passing through $0$ precisely like $-x/y - Y/X$. We assume in
the discussion that 
$Y/X$ is not affected by the movement (Lie's point transform). Note also, that if
$Y/X=\rho>1$, we have if $U \rightarrow I$, that $\rho>1$ (hyperbolic) for degenerate points,
why if we limit ourselves to elliptic approximations, we do not have degenerate points.

The invariance principle for movements, has a correspondent principle in operator space. 
Note that if we start from the parametrices as 
Fredholm operators, with a decomposition $N(E) \bigoplus D(E)$, we see that modulo $C^{\infty}$, 
that parametrices to hypoelliptic operators have $N(E)=\{ 0 \}$. Orthogonality for the symbol space 
induces a corresponding relation in the operator space. If for two analytic symbols $f_{1},f_{2}$, we assume 
$f_{1} \prec f_{2}$, we note that a necessary condition for inclusion for the correspondent space of 
operators is $f_{1} \prec \prec f_{2}$. Thus, we have existence of $N$, such that for the operator space 
$ \frac{d}{d x} f_{1}^{N} \prec f_{2}^{N}$. In particular, we can write $\frac{d}{d x} \log f_{1}^{N} \rightarrow 0$ 
for some $N$ and simultaneously $f_{1}^{N} \prec f_{2}^{N}$. Note that for a polynomial, we have 
always $\frac{d}{d x} P \prec P$. Further, if $\mbox{ rad }f_{2}  \subset \mbox{ rad }f_{1}$,
we have $N_{1} \subset N_{2}$, for the corresponding zero space. Radon Nikodym's theorem can be used, 
$I(f_{1})=I(g f_{2})$, for $g \in L^{1}$. The conclusion is that a necessary condition for inclusion 
of the corresponding operator spaces, is that one symbol ideal is strictly weaker that the other and 
a strict dominance of symbol ideals implies an inclusion of operator spaces.
It is for this sufficient to consider the phase space, $\log f_{1}^{N} \prec \log f_{2}^{N}$ 
and $\frac{d}{d x} \log f_{1} \rightarrow 0$ in the $\infty$ implies $(I_{2}) \subset (I_{1})$. 
Note that (\cite{Schwartz66}) $\mathcal{D}$ is dense in $\mathcal{D}_{L^{p}}$ and in 
$\dot{B}$, but not in $B$. We have that $(\dot{B})'=\mathcal{D}_{L^{1}}'$ and $(B)'=\mathcal{D}_{L^{\infty}}'$. 
Further we have that $(\dot{B})'$ is the limit in $B'$ of $\mathcal{E}'$.

Concerning $f \in (I_{1}) \subset (I_{He})^{\bot}$, we assume $(I_{He})$ with a global pseudo base
and $< f, d \mu>=0$ continued to $< \tilde{f}, d \mu >=0$. Assume the 
continuation $\tilde{f} \sim e^{g} f$ with $g \in L^{1}$, that is $f \bot e^{-g} d \mu$. If $dv$ 
corresponds to $(I_{Phe})$, such that $dv \sim e^{-g} d \mu$, then $e^{g} d v$ is downward bounded. 
 Another example is given by $d \mu$ of type $0$, 
$e^{- g_{1}} d \mu$ of type $A_{1}$, further $e^{-g} d \mu$ of finite type. A final example 
is given by $d \mu$ with a trivial kernel.

Assume $\gamma_{1}=\gamma_{2}$, on $\Delta$ and construct a neighbourhood $J U \gamma_{1}=V J \gamma_{1}=V \gamma_{2}$.
We have $U \gamma_{1}^{2}=(U \gamma_{1})^{2}$, but $V \gamma_{2}^{2} \neq (V \gamma_{2})^{2}$ and
$V (e^{v} \gamma_{2}) \sim J U \gamma_{1}$ and $J \gamma_{1}^{(2)}=e^{v_{2}} \gamma_{2}^{(2)}$, 
we then have $v_{2} \neq v_{1}$, but $v_{1}=v_{2}$ on $\Delta$. For instance: 
$\{ 0 \} = \Delta (\gamma_{2}^{(2)}) \subset \Delta(\gamma_{2})= \Delta(\gamma_{1})=\Delta(\gamma_{1}^{(2)})$.
Thus, $\{ e^{v_{2}} - 1 \} \subset \{ e^{v_{1}} - 1 \}$, that is $\{ v_{2}=0 \} \subset \{ v_{1}=0 \}$
and for the correspondent geometric ideal, we have $I_{1} \subset I_{2}$. When we assume $U$ linear 
in $\gamma$, we do not assume $U$ simultaneously linear in $\zeta$, for $V$ we do not assume $V$ linear in 
$\phi$ or $\zeta$. $V = J U J^{-1}$ and $I = I_{Hyp}$ implies $V I^{\bot}=J U I$.

The boundary is defined by first surfaces $S$, invariant for all movements 
and $S_{0}=S \backslash \{ x_{0} \}$. 
Existence of regular approximations is guaranteed 
by $\Delta$. Transversals are associated to reflection axles  in $S$, we consider one movement at a time. 
$\Delta$ is represented by $U^{\bot}=I$ and the choice of movement determines the properties of the neighbourhood of $\Delta$. 
When $\gamma$ preserves a constant value in the $\infty$, we have $\gamma(x,y) \sim P(\frac{1}{x},\frac{1}{y})$, 
as $x,y \rightarrow \infty$, for a polynomial $P$ and we can determine $\bot$ as independent of $\mid x \mid,\mid y \mid \rightarrow \infty$. 
In the case when the regular approximation does not have a reduced $\bot$ measure, that is 
$\int g_{reg} d \mu=0$ where $d \mu$ is not reduced, we must take into account orbits among the
possible approximations. 
Sufficient for this to occur, is that we do an adjustment with point support measure  
in Cousin's continuation over the boundary (\cite{Dahn13}, cfr the last section in this article). Thus, for a reduced representation in 
the moment problem for the measure, we can assume transversals without orbits or orientation of orbits.

Assume $f=e^{\phi}$, and consider the problem when $U$ preserves analyticity. If $\big[ U,I \big]=\big[ I,U \big]$ 
and $U$ is acting linearly in the phase, we have that $U(f_{1}f_{2})=e^{U \phi_{1} + U \phi_{2}}
=U(f_{1}) U(f_{2})$ (``point topology''). Note that if $U=I$ on $W$ planar in $\mathcal{O}_{AD}$ (\cite{AhlforsSario60}),
we have that $U$ linear. 
When $U=I$ on $W \ni \infty$ planar,  we must assume $(X,Y)_{T} \rightarrow (X,Y)_{1/T}$ continuous. 
When $U \rightarrow V$, it is not sufficient to consider tangents and we assume $F_{T} \rightarrow F_{1/T}$ 
preserves continuity over the axes for invariance.

Note that when $U$ is linear over $\phi$, we can assume $U(\overline{\phi})=\overline{U(\phi)}$,
$U(i \phi)=i U(\phi)$ and $U(- \phi)=-U(\phi)$. In particular this can be assumed when $\phi$ defines a planar domain 
$\phi_{1} \bot \phi_{2}$ in $\mathcal{O}_{AD}$ (or standard complexified). For $J \phi$, we have 
for some iterate $\phi_{1}^{N} \bot \phi_{2}^{N}$, however we also have a $e^{v}$ for $v \in L^{1}$ 
according to Radon-Nikodym's theorem, that is we have a non planar domain. 
Since $v$ is determined by the movement, we will argue that it is sufficient to consider a domain for $v$
on one side of a hyperplane.

\newtheorem*{res1}{Invariance principle}
\begin{res1} 
 By considering movements as functionals, we can uniquely relate movements on the hyperboloid to movements
 in an euclidean metric. We assume $J : \gamma_{1} \rightarrow \gamma_{2}$ continuous and $\gamma_{1}=\gamma_{2}$
 on a set $\Delta \neq \{ 0 \}$. If we assume $\gamma_{j}$ analytic and $\gamma_{2}$ reduced
 with respect to $\gamma_{1}$, we can represent $F(J \gamma_{1})=F(e^{v} \gamma_{2})$ and $v \in L^{1}$
 The boundary measure represents a very regular distribution $B$ in $\mathcal{D}_{L^{1}}'$, with kernel 
 $\widehat{B} \subset L^{1}$.
 \end{res1}

That is we do not assume $\gamma_{j} \in L^{1}$, but $\frac{d \gamma_{2}}{d \gamma_{1}} \in L^{1}$ with
respect to the boundary measure. Thus, since $\gamma_{2}$ is reduced, $\frac{d v}{d \gamma_{1}}$ is in $L^{1}$
and we can conclude that $v \in L^{1}$ with respect to the boundary measure. The boundary distribution is
thus of real type (\cite{Martineau}) and it is sufficient to consider $\mathbf{R}^{n}$ for the definition of boundary 
condition. If we assume $B(e^{v})=\widehat{B}(v) \in L^{1}$
implies $v \in L^{1}$, we have that $v \in L^{1}(\mathbf{R}^{n})$.
Concerning the boundary condition, we note that the condition on a very regular boundary $\delta_{0} - C^{\infty}$ 
does not imply analyticity. The proposition that $\mbox{ ker }E$ (parametrix) can be replaced by 
$C^{\infty}$, is implied by homogeneous hypoellipticity.

\section{Movements}

Assume that $\pi$ is a plane through the origo and $K$ is the light cone. We define 
$V_{R} : \quad \pi \cap K=\{ 0 \}$,$V_{L} : \quad \pi \cap K=\{ L \}$,$V_{T} : \quad \pi \cap K=\{ L_{1},L_{2} \}$.
Assume $R$ reflection points (invariant points). We then have $R \in V_{T}$ implies elliptic 
movements, $R \in V_{L}$ implies parabolic movements, $R \in V_{R}$ implies hyperbolic movements.

Every movement on the hyperboloid, can be given as reflections with respect to a plane through a fixed point.
They are divided into direct movements: reflection with respect to a plane through an axes and indirect 
movements: reflection with respect to a plane through a point. We can for instance assume the axes 
$x=y$ and the point is 0.

The lineality $\Delta$ is defined in the domain $\Omega$ through translation invariance for an analytic symbol. 
The corresponding movement in lorentz geometry is $I$, that is 
$I \gamma=\gamma$. More precisely $\tau F(\gamma)(\zeta)=F(\gamma)(\zeta)$
or equivalently $F(U \gamma)=F(\gamma)=F(I \gamma)$. Assume $U \rightarrow I$ is a movement and 
let $U \rightarrow \tilde{U}$ be the mapping $(I) \rightarrow \Omega$. Given $U \gamma$ 
analytical, we have that $\tilde{U} \zeta$ defines a neighbourhood of $\Delta$.
In the same manner, if we let $U^{\bot} \rightarrow I$,
with $U^{\bot} \rightarrow \tilde{U}^{\bot}$, this does not imply $U^{\bot} \gamma \in H_{m}$, but given
$U^{\bot} \gamma$ analytic, we have that $\tilde{U}^{\bot} \zeta$ is continuous. Note that without the  condition 
on analyticity, we do not have $\lim_{U \rightarrow I} \tilde{U} = \lim_{U^{\bot} \rightarrow I} \tilde{U}^{\bot}$
For this reason, we define $\Delta=\{ \zeta \quad I \gamma = \gamma \quad \gamma \in (I) \}$, where 
$(I)=\{ \gamma \quad F(\gamma)(\zeta)=f(\zeta) \}$ for some $F$. $\Omega$ can be defined through 
$\{ \zeta \quad U \gamma(\zeta) \mbox{ close to } \gamma(\zeta) \quad \zeta \in \Delta \}$.
In the same manner, we can define $\Delta^{\bot}=\{ \zeta \quad I \gamma=\lim_{U^{\bot} \rightarrow I} U^{\bot} \gamma \}$.
Thus, we can assume $U \gamma \subset \{ F(\gamma)=const \}=S$ first surfaces with compatibility conditions 
$U^{\bot} \gamma \cap S_{0} = \emptyset$ where $S_{0}=S \backslash \{ \zeta_{0} \}$ and $\zeta_{0}$ 
a singular point, further that $U^{\bot} \gamma$ is analytic close to $S$. 
The movement $U_{e}$ (elliptic), has $U^{\bot}_{e}$ as ``transversals'', $U_{p}$ (parabolic) has
$U^{\bot}_{p}$ of the same character as $U_{p}$. For $U_{h}^{\bot}$ (hyperbolic) we have 
$U_{h} \gamma \subset S$ and $U_{h}^{\bot} \gamma$ approximates $S$. 
Further by the Fourier dual $(U_{h})^{*} \simeq (U_{e})^{\bot}$. Correspondingly in euclidean metric 
we have $(V_{\tau})^{*} \simeq (V_{e})^{\bot}$ where $\tau$ is translation and $e$ is rotation. 
We assume here $(I)=(I_{Hyp}) \oplus (I_{Phe})$ are
analytic functions, that is we assume $(I_{Hyp}) \ni \gamma \rightarrow J \gamma \in (I_{Phe})$ 
preserves analyticity.
Thus, we have $U_{e} \rightarrow U_{h}^{\bot}$
where transversal lines can be traced by translation, $U_{e}^{\bot} \rightarrow U_{h}$ first surfaces
(multivalued) can be traced through planar movements. 

However we do not have that $U \gamma \in (I_{Hyp})$ implies $U^{\bot} \gamma \in (I_{Hyp})$, 
further we do not assume
$U^{\bot} \gamma \subset (I_{Hyp})^{\bot}$ or $V^{\bot} J \gamma \notin (I_{Phe})$. We define $J$ 
through $L(\gamma, U^{\bot} \gamma)=0$ implies $E(J \gamma, V^{\bot} J \gamma)=0$. Thus, $L(U \gamma,U^{\bot} \gamma) \equiv 0$ 
for every $\gamma$. We use $<,>$ to denote respective scalar products.
For a regular (reversible) movement (axes $\notin H_{m}$) we have that $\big[ < \big[ I,U \big] \gamma,\xi>=0 \Rightarrow \xi=0 \big]$
implies $U=I$. In particular, if $< (U - U^{\bot}) \gamma,\xi>=0$ implies $\xi=0$ then $U=U^{\bot}$ 
that is we have points in common. Through the compatibility condition, we can assume points in common 
on $S$. Note further $\Delta$ is joint for $\gamma,J \gamma$.

Given that the joint points can be defined through $\lim U^{\bot} \gamma = \lim U \gamma$ (with orientation)
we see that ${}^{t} U U^{\bot} - I \equiv 0$ in limes. Given $U \gamma \in \Sigma=\{ F(U \gamma)=const \}$,
by reverting the orientation for $U^{\bot} \gamma$ the continuation $U \gamma - U^{\bot} \gamma$ 
can be taken in the sense of Cousin (\cite{Dahn13}). We can assume product topology for $U \gamma$, 
the continuation using $J$ is continuous and in $L^{1}$ (with respect to the boundary measure).
If ${}^{t} U \gamma \neq U^{\bot} \gamma$ except for a discrete set we can regard $U^{\bot} \gamma$
as a continuation of $U \gamma$. Note that we have existence of $\gamma$ with $U \gamma = I \gamma$ 
implies existence of $\gamma'$ such that $U^{\bot} \gamma'=\gamma'$ with $\gamma \neq \gamma'$. 
Thus $U^{\bot} U \gamma = I U \gamma$
and $U U^{\bot} \gamma ' = U I \gamma'$ that is on invariant points we have $\big[ U,I \big]=\big[ I,U \big]$
(for instance $U \gamma=\gamma'$).

Further $< U \gamma, U^{\bot} \gamma> \equiv 0$ iff $< \gamma, {}^{t} U^{\bot} U \gamma> \equiv 0$
why ${}^{t} U^{\bot} U \in {}^{\circ}(I)$ (annihilators). The movements are primarily considered 
in $H'$, though using the moment problem, if the movement is analytic in a set $E_{0}$, it can be
continued to $C$, assuming the compatibility conditions above.

Assume orthogonality is defined by
$<f,g> = I_{f}(g)=0$. We then have that $N(I_{f})$ is defined for $g \in H$, if $f$ has regular 
kernel. If $N(I \oplus I^{\bot})=\{ 0 \}=N(I) \cap N(I^{\bot})$, we can write 
$V \cup V^{\bot}=\Omega$ and $V \cap V^{\bot}=\{ 0 \}$. 

When two mirrors are used,  we get a non-commutative group. Consider the reflection 
$S=(S_{1},S_{2})$ through the diagonal
$\pm x=y$, $R$ through the real axes and $T$ through $0$. Note that the diagonal has two 
generators $S_{1},S_{2}$ and
$R S_{1} \neq S_{1} R$. But we have that $R S_{1} = S_{2} R$ and $TR S_{1}= R S_{1} T$. 
Thus ${}^{t} T=T$,${}^{t} R=R$ and
${}^{t} S_{1} = S_{2}$. For $z \rightarrow z^{*}$ we have that $TR S_{1}=RT S_{1} = R S_{1} T = T S_{2} R = S_{2} T R =
S_{2} R T$. That is $z^{*}=-i z$. Note that for harmonic conjugation that for a closed form 
$\varphi$, if $\varphi^{*}=-i\varphi$ that $\varphi=\alpha d z$, where $\alpha$ is analytic 
locally.

A joint boundary for $(I_{1}),(I_{2})$ is represented using first surfaces 
$\Sigma=$ $\{ F(U \gamma_{1})(\zeta)=$ $F(\gamma_{1})(\zeta) \}$ for every movement $U$ and 
$\Sigma'=\{ F(V \gamma_{2})(\zeta)=$ $F(\gamma_{2})(\zeta) \}$, 
for every movement $V$ and $\Sigma \simeq \Sigma'$. Thus, where 
$J_{\delta} : \gamma_{1} \rightarrow \gamma_{2}$ and $J_{\delta} U=V J_{\delta}$
we assume $F(J_{\delta} U \gamma)=const.$ for all $\delta \rightarrow 0$. Consider the compatibility 
condition $F(V^{\bot} \gamma_{2}) \neq 0$, for some $V^{\bot}$ close to and at distance from the boundary. 
$V^{\bot} \gamma_{2}=J U^{\bot} \gamma_{1}$ is defined so that $F(V^{\bot} \gamma_{2}) \rightarrow \Sigma'$ 
iff $F (J U^{\bot} \gamma_{1}) \rightarrow \Sigma'$ regular. This is interpreted as two-sided, 
that is the invariance principle in this case is extended with a compatibility condition. 

\newtheorem*{res2}{Involutive movements}
\begin{res2}
 Any movement on the hyperboloid is involutive. In our mixed model, any involutive movement in euclidean
 metrics, has a correspondent movement on the hyperboloid. Orthogonal movements, are not necessarily
 involutive, but our invariance principle maps $U^{\bot} \rightarrow V^{\bot}$ uniquely.
\end{res2}

\subsection{The lifting operator}

Consider $L(Ux,y) \rightarrow E(V x.y)$, where $E$ denotes the euclidean metric, through 
the invariance principle. Define $\{ V x \}^{\bot}=\{ y=Rx \quad y \bot V x \}$ 
or $< V x, y >=R(V x)$, where $R$ is defined using Radon-Nikodym's theorem 
and is specific for the movement. We define $V^{o}=\{ R \quad R(Vx) = 0 \quad \}$, a closed set 
that is generated by $V$. Conversely, if $L({}^{t} R_{U} U x,x)=0$ $\forall x$, which implies 
$E({}^{t} R_{V} V x,x)=0$ $\forall x$, where $R_{U}U=U^{\bot}$ and $R_{U}: X' \rightarrow X'$. 
Thus, $\big[ L_{x}, U^{\bot} \big] \sim \big[ E_{x}, V^{\bot} \big]$
and over the reflection axes 
$J(U^{\bot})=(JU)^{\bot}$.

Using the condition that $f$ and $g$ have lineality in common, we see that we have a continuous mapping
$\{ f-c \} \rightarrow \{ g-c \}$, through a planar (transversal) mapping.
In particular, with the compatibility conditions, we have $sng f \rightarrow sng g$. Note that we assume $\log f,\log g \in L^{1}$.
A reflection axes is defined by $R_{1} : {}^{t} r \gamma_{1}=\gamma_{1}$. Given $J : \gamma_{1} \rightarrow \gamma_{2}$
and $J {}^{t} r_{1}= r_{2} J$, there is a correspondent set for invariance $r_{2} \gamma_{2}=\gamma_{2}$
that is defined by $R_{2}$. For $\gamma_{1}$, we assume symmetry with respect to $R_{1}$. The corresponding 
proposition for $\gamma_{2}$, is dependent on $e^{v}$ and symmetry for $v$ with respect to $R_{2}$.

Using Radon-Nikodym's theorem $({}^{t} r_{1} -1) \gamma_{1}=0$ and $({}^{t} r_{2} -1) \gamma_{2}=0$,
we then have $({}^{t} r_{2} -1) J e^{v}= ({}^{t} r_{1} -1)$. Assume $y=y(x)$ and consider 
$k(x,y)=\frac{y}{x}$. 
Given $k > 1$ and $\frac{d y}{d x} - \frac{y}{x}>0$, this implies $1 < \frac{y}{x} < \frac{d y}{d x}$ (elliptic movement). 
Given $k < 1$ and $\frac{d y}{d x} - \frac{y}{x} < 0$, this implies $\frac{d y}{d x} < \frac{y}{x} < 1$ 
(hyperbolic movement). 
Assume $r_{2} d \gamma_{2} \sim d r_{2} \gamma_{2}$, we then have 
$\frac{r_{2} d \gamma_{2}}{r d \gamma_{1}}=e^{v}$ or $\frac{r_{2} (X_{2},Y_{2})}{r (X_{1},Y_{1})} \sim e^{v}$
where $\frac{d \gamma_{2}}{d t}=(X_{2},Y_{2})$. Given $\frac{k(X_{2},Y_{2})}{k(X_{1},Y_{1})}=const$
and $X_{1} < Y_{1}$, we have that $X_{2} < Y_{2}$ and so on. Note that $\gamma_{2}$ can be seen 
as a continuation of $\gamma_{1} \subset (I_{He})^{\bot}$. Given $J r \sim {}^{t} r_{2} e^{v}$ 
and $< r \gamma_{1},\gamma_{2}> = < \gamma_{1}, r_{2} \gamma_{2}>$ and $< r \gamma_{1}, \gamma_{2} > = < {}^{t} r_{2} \gamma_{1},\gamma_{2}>$
that is $r \sim {}^{t} r_{2}$. Assume $A$ an annihilator for $(I_{1})$ and $B$ an annihilator 
for $(I_{2})$, both closed operators. For instance $A=r-1$ and $B={}^{t} r_{2} -1 $. We then have 
$(I)=(\mbox{ ker }A)$ and $(J)=(\mbox{ ker }B)$. We have that if $N(I)=N(J)$, given an analytic representation 
of $(I),(J)$, that $\mbox{ rad }(I) \sim \mbox{ rad }(J)$. For movements, if we have analytic representations of 
both $r, {}^{t} r_{2}$ and $r-1,{}^{t} r_{2} -1$, then the geometric sets coincide.
 
\newtheorem*{res3}{Fundamental representation}
\begin{res3}
 We assume given $f$ analytic that $\gamma_{j}$ are analytic and that the equation $F(\gamma_{j})(\zeta) \sim f(\zeta)$
 can be solved in $\mathcal{D}_{L^{1}}'$. When we can restrict the model to polynomial $\gamma_{j}$,
 $F$ can be constructed from $f$ and the parametrices to $\gamma_{j}$.
 In particular, where we have isolated singularities, we can represent $F$ as a measure with compact support.
 \end{res3}

A very regular boundary is characterized by (\cite{Parreau51}) singularities located in a locally finite set of
isolated points or segments of analytic curves. Particularly, when the regular approximations have isolated
singularities for some higher (finite) order derivative. When the model is considered in $\mathcal{D}_{L^{1}}'$,
the real Fourier transform can be represented $P(\zeta) \tilde{f}_{0}(\zeta)$ (\cite{Dahn13}), for a polynomial $P$
and $\tilde{f}_{0}$ very regular. The fundamental representation can be derived in several ways. The parametrices
to hypoelliptic differential operators are very regular in the sense that the Schwartz kernel is regular outside the 
diagonal. Assume for instance $E$ a parametrix to $\gamma(D) \delta$, where $\gamma(\zeta)$ is a hypoelliptic polynomial.
Then $\big[ E, \gamma(D) I \big] \sim \gamma(D) \big[ E, I \big]$ has Fourier-Borel transform $F(\zeta)\gamma(\zeta)$,
where $F$ is $\widehat{\big[ E,I \big]}$ and $\big[ E,I \big]=\big[ I,E \big]$, in this setting, when the parametrix is two-sided. 
Note also that $F(\gamma(D) \varphi)=E(\gamma(\zeta) \widehat{\varphi})$, when $\widehat{\varphi} \rightarrow 1$.
\subsection{Reflections}

All movements in $(I_{Phe})^{\bot}$ can be represented as reflections. Singularities for 
$(I_{Phe})$ can all be related to the boundary. Every point on the boundary can be reached through 
reflections in $(I_{Phe})^{\bot}$ emanating from an ``origo'' on the boundary. If the points are 
reached through planar reflection, then analyticity is preserved.
We can represent points that can be reached through hyperbolic movements (translation), parabolic movements (scaling)
elliptic movements (rotation) in the same leaf. 
In cylindrical domains (order 0), the translations are parallel with one axes.
Pseudo convex domains are locally cylindrical at the boundary. 

Given $P,Q \in  H_{m}$ (the web of the hyperboloid), we can prove existence of unique movement 
$P \rightarrow Q$, that corresponds to a unique (euclidean) translation $P \rightarrow Q$. 
If $\pi$ is a plane through $P,Q$ and $R$ a line $\bot \pi$ implies a unique existence of a 
rotation axes $R \sim PQ$ (line). $R$ is not uniquely given by the movement, 
we have $\infty^{1}$ possibilities in $V_{T}$, $\infty^{1}$ in $V_{R}$ and $\infty^{2}$ in 
$V_{L}$ (\cite{Riesz43}). However, given that $R$ is a hyperbolic movement, 
it has a unique representation as a (euclidean) translation determined by $P,Q$.

\subsection{Distance functions}

Concerning the invariance principle, assume $d_{\Gamma}=d_{x} \otimes d_{y}$, where $d_{\Gamma}$ denotes
the distance to $x=y$ in $\gamma_{1}$. We now have three cases $d_{x}/d_{y}=const$ ($d_{\Gamma} = const$)
$d_{x}<d_{y}$ and $d_{x}>d_{y}$. The corresponding distance function 
in $\gamma_{2}$ can be induced according to
$(J^{-1})^{*}d_{\Gamma}$ ( pseudo distance.)

The normal is defined in $(I_{1})$ starting with the tangent. 
In the case where $(I_{2})$, we have blow-up according to
$df \rightarrow f$, why we prefer to start with invariant points. 
On first surfaces where all points are invariant, we consider transversals 
on the form $V^{\bot} \gamma$.

Consider $(I_{1}) \subset (I_{2})^{\bot}$ and $d=(v_{1},v_{2})$, that is product topology. 
We assume Runge's property for $d( \gamma_{1} - \gamma_{2}) \sim \inf (v_{1}(\gamma_{1}),v_{2}(e^{v} J \gamma_{1}))$ 
for $v \in L^{1}$. Assume for instance $e^{v} \gamma_{2}=z e^{v_{1}(z)} \sim e^{w(\frac{1}{z})}$ 
and we have existence of $\eta$, such that $w(\frac{1}{z}) \sim \frac{1}{\eta(z)}$, where 
$\eta \in L^{1}$. When $f = e^{\phi}$, we obviously have $\big[ I, e^{v} \big](e^{\phi})=I(e^{\phi + v})$ 
$ = \big[ e^{v},I \big](e^{\phi})$.

 Assume that distances are given by $(v_{1},v_{2})$, we then have that for invariant points that 
$v_{2} \leq v_{1}$ for a hyperbolic movement, $v_{1} \leq v_{2}$ for an elliptic movement. Assume 
$v_{1} \mid f \mid \leq 1/v_{2}$, where $v_{2}$ is algebraic and $v_{1}=const$, we then have that 
$v_{2} < 1$ for $\widehat{f} \in C^{N}$ and $v_{2} > 1$ for $\widehat{f} \in \mathcal{D}_{L^{1}}'$ 
(finite order). As $v_{1} \neq const.$ we have $v_{2}^{2} \mid f \mid \leq v_{1}v_{2} \mid f \mid \leq C$, 
when $v_{2} \leq v_{1}$ and $v_{1}^{2} \mid f \mid \leq v_{1} v_{2} \mid f \mid \leq C$ in the converse 
case. When $\deg v_{2}^{2} = N$, we have $\widehat{f} \in C^{N}$. For instance, $v_{1}(x,y)=const$ 
we have $v_{x1}(z)=v_{x1}(\overline{z})$ and $v_{x2}(y-x)=v_{x2}(x-y)$, that is we assume $x$ fixed 
and $v_{1}v_{2}=v_{2}v_{1}$. 

We define $v$ through movements $J U \gamma_{1}=V \gamma_{2}$. Assume $(d_{1},d_{2})$ 
distances to the euclidean axes, that is $d_{1} < d_{2}$ on the reflection axes, implies $V$ 
translation, $d_{2} < d_{1}$ implies $V$ rotation and $d_{1}=d_{2}$ implies $V$ scaling. 
Assume $E_{0}$ with product topology and $d \mu \bot E_{0}$. The Runge property implies that 
$d \mu \bot C$. Thus, if we have separate invariance in $x$ and $y$, we have
invariance in $(x,y)$.

Concerning monotropy, we have a separate (pluri complex) condition that is represented by $\epsilon-$ 
translation. For instance if $v = \big[ v \otimes I + I \otimes v \big]$, it is sufficient that 
$\big[ I,e^{v} \big](f)=\int \big[ I,e^{v} \big](x,y) f(y) d y=\int ( \int I(x,z) e^{v}(z,y) d z) f(y) d y=\int e^{v}(x,y) f(y) d y$.
Assume for $x$ fixed, $G(f)(x)=\int e^{v(x,y)} f(y) d y$ and for $y$ fixed, 
${}^{t} G(f)(y)=\int e^{v(x,y)} f(x) d x$.
Thus $< G(\gamma_{1}), \phi>=< \gamma_{1}, {}^{t}G(\phi) >$. Given 
$J U \gamma_{1}=V \gamma_{2}$
and $J \gamma_{1} = \gamma_{2}$, which implies $U \gamma_{1}=\gamma_{1}$ and this implies $V \gamma_{2}=\gamma_{2}$.
Given $\gamma_{1} \bot J^{-1} \gamma_{2}$ with $F(J \gamma_{1})=F(e^{v} \gamma_{2})$, thus if 
$U \gamma_{1} \bot V \gamma_{2}$, we must have $v=const$ on all invariant sets.

If we consider $(\gamma_{1},\gamma_{2}) \in (I_{Hyp}) \times (I_{Phe})$, 
we can consider an associated distance $d=d_{1} + i d_{2}$. 
In the case with a regular complement and $\gamma_{1} \rightarrow \gamma_{2}$ with joint lineality, 
if $d_{2,N}$ denotes the distance to the lineality for $\gamma_{2}^{N}$, we have $d_{2,N} \geq d_{1,N}$. 
Thus in the case when $\Omega \backslash \Delta_{N} \downarrow \{ 0 \}$, we have that $d_{2,N} \uparrow$, 
as $N \uparrow$. 

Alternatively we can give  
$(x , y^{N}) \in (I_{Hyp}) \times (I_{He})$. We assume $\frac{1}{d_{2}}$ a pseudo distance
and also a distance. 
Given the condition $d_{1}/d_{2} \rightarrow 0$, 
when $y \rightarrow 0$, we have that $\frac{d_{1}d_{2}}{d_{1}^{2} - d_{2}^{2}}=\frac{1}{\frac{d_{1}}{d_{2}} - \frac{d_{2}}{d_{1}}} \rightarrow 0$, 
as $y \rightarrow 0$. Thus if $\tilde{d}^{2}=(d_{1} + i d_{2})^{2}$, we have that 
$\frac{\mbox{ Im } \tilde{d}^{2}}{\mbox{ Re } \tilde{d}^{2}} \rightarrow 0$, 
as $y \rightarrow 0$. 

Assume  $\pi$ a line (axes for reflection) $= \{ d_{1}(z)=0 \}$ and $d_{1}(t z)=t d_{1}(z)$, 
$t>0$. 
If $\pi$ is a plane, we can use a two-sided distance 
function, for instance $\{ d_{1}=const \} \cup \{ d_{2}=const \}$, that corresponds to 
$\overline{N(d_{1})} \cup \overline{N(d_{2})} \sim \overline{N(d_{1} d_{2})}$.
Concerning first surfaces to distance functions, given $d^{2}$ algebraic (limit of algebraic functions),
when $d^{2}$ is a distance, it is locally 1-1. Given Schwartz type topology, we can assume these first 
surfaces have the same properties as the zero sets.

Assume $S_{1}=\{ \zeta \quad d_{1}=0 \}$ and $\Sigma_{1}=\{ \gamma_{1} \quad F(\gamma_{1})=const \}$, 
then $\gamma_{1}$ is locally bounded, we have thus $d_{\Sigma} \leq c d_{\Gamma}(f)$, that is 
$d_{\Gamma}=0$ implies $d_{\Sigma}=0$. Assume $d_{1,\Sigma} \rightarrow d_{2,\Sigma}$. 
If we assume $J_{\delta}$ continuous, we have $d_{\Gamma}=0$ implies $d_{2,\Sigma}=0$. 
conversely, if $J_{\delta}^{-1}$ is continuous, we do not necessarily have $d_{2,\Gamma}=0$.

\subsection{Singularities}
$\Delta$ is defined through translation in $\Omega$, which corresponds to $\Delta=\{ I \gamma=\gamma \}$.
We define a neighbourhood of $\Delta$ with respect to a lorentz movement $U \rightarrow I$ 
and $U \rightarrow \tilde{U} \zeta \in \Omega$. In the case where $U \gamma \rightarrow I \gamma_{0}$, 
then we can regard a singularity as an isolated point on $\Delta$, given that $U$ preserves 
analyticity (planar), the point is reached using $\Delta$ as a strict carrier for the limit, 
that is it is of no consequence what neighbourhood we consider. We assume all through this article 
that $f(\zeta)=F(\gamma)(\zeta)$ and $\log \mid f \mid \in L^{1}$, that is the singularities
are given by $\log \mid f \mid$ and are of finite order.

In the case with $U_{e}$ elliptic, we have that $\mid V_{e} f \mid=\mid f \mid$, that is 
these changes of variables do not affect the singularities. In this case, the singularities on $\mid f \mid=1$ 
can be seen as isolated and are completely defined by transversals. In the case where $\mbox{ Re }f \bot \mbox{ Im }f$,
we can start with an arbitrary point on $\mid f \mid = 1$, and approximate the singularity 
by elliptic rotation, that is $U_{e}$ such that $\mid V_{e} f \mid=1$ and $\tilde{V}_{e} \zeta = \zeta_{0}$. 
Given $\mid f \mid^{2}$ analytic (does not imply $f$ analytic), we can assume $\Omega=\{ \zeta \quad \mid f \mid=1 \}$ 
connected ($\mbox{ Im }f \sim \mbox{ Re }f$). 
Consider the movement as a functional and $V_{e} F(\gamma)=F({}^{t} V_{e} \gamma) \sim V_{e} f$. 
Then $< V_{e} \gamma,\gamma>=1$ implies $< \gamma, {}^{t} V_{e} \gamma>=1$, that is we can 
assume ${}^{t} V_{e}$ rotation (given a normal model). 

Concerning parabolic movements $U_{p} f=F({}^{t} U_{p} \gamma)$ and $\frac{\widetilde{y}}{\widetilde{x}}=\frac{y}{x}=\rho$
and $y=y(x)$, ${}^{t}U_{p} y=\widetilde{y}$. We are considering $\Phi(\frac{y}{x})=(t x, t y)$ for $t$ real. 
It is for this reason sufficient to consider singularities that are given by $\eta(x)=y(x)/x$, 
where $\eta$ is algebraic iff $y$ algebraic. Note that $\{ \eta=\frac{d \eta}{d x}=0 \} \subset \{ y=x \frac{d \eta}{d x}=0 \}$
We have that $x \frac{d \eta}{d x}=0$ iff $\frac{d y}{d x} \sim -x/y$ (\cite{Dahn13}) 
Note if $(x,y) \rightarrow (x',y')$ with $\frac{dy}{dx}=\frac{dy'}{dx'}$ and $y'(x)=y( \alpha x)={}^{t}A y (x)$,
then if ${}^{t}A=A$, $ \frac{d {}^{t} A y}{d x}=A \frac{d y}{d x}$ implies assuming ${}^{t}A =A$, 
$\frac{d {}^{t }A}{d x}=0$.

For translation $V_{h} f=F({}^{t} V_{h} \gamma)$, in the moment problem, 
assume that $A$ defines regular approximations and $V_{h} \gamma$ 
defines a movement on the first surface to $f$. Then,
$\int A d (V_{h} \gamma)=\int A d \gamma$, that is $A$ is not dependent on $V_{h}$. 
In this case the singularity is not 
affected by $V_{h}$. If $F$ is linear and $F : const \rightarrow const,$, we can assume 
$\{ F(V_{h} \gamma)=0 \} \sim \{ F(\gamma) - \lambda' \}$,
that is first surfaces to analytic functions, with regularity conditions as with (\cite{Nishino68}). 
In particular, when $\lambda \rightarrow \lambda_{0}$, there are $V_{h}$ such that 
$V_{h} \gamma \rightarrow \gamma_{0}$. 
For $\tilde{V}_{h} \zeta$ we can compare with Abel's problem (\cite{Julia24}). 
Given isolated singularities, we can assume $\tilde{V}_{h} \zeta=\zeta + \eta$, 
for some $\eta$.

Concerning the two mirror model, we consider $A \rightarrow \gamma_{1} \rightarrow \gamma_{2} \rightarrow B$
where $A,B$ are situated on first surfaces to $f$. In the planar case, where $\mid \gamma_{1} \mid \leq 1$,
$\mid \gamma_{2} \mid \leq 1$ and $\gamma_{1} \bot \gamma_{2}$, then every  $B$ can be reached, 
independent of starting point $A$ on a first surface. Thus in the planar case, 
every pair of points on the first surfaces, can be combined using a continuous path.  

Consider $U \rightarrow V \in (I_{2})'$ and $(I_{1}) \subset (I_{2})^{\bot}$, that gives 
a continuation $\phi$ of $\gamma_{1}$, that is $< \phi,\gamma_{2}>=0$. Further 
$< V \phi, \gamma_{2} >=< \phi,{}^{t} V \gamma_{2}>$ defines
${}^{t} V \in (I_{2})'$ and ${}^{t} V \gamma_{2}=\gamma_{2}$ iff 
$V \phi = \phi$ 
and the invariance principle is here considered in $\mathcal{D}_{L^{1}}'$. Assume 
$J : U \rightarrow V$ and $F(J U \gamma)={}^{t} J F(\gamma)$. 
A chain given by $U$ is mapped by $J$ onto a chain given by $V$.

Concerning algebraicity, assume $\mid f \mid < C \mid P \mid$ in $\infty$, where $P$ is polynomial, 
then for $P$ to serve as a weight in $\mathcal{D}_{P}'$, it is necessary that $I \prec \prec P$. 
We then have that $P_{1} \prec \prec P_{2}$ implies $\mathcal{D}_{P_{2}}' \subset \mathcal{D}_{P_{1}}'$ 
which corresponds to analytic functionals of finite type. In particular, for a reduced polynomial, we have 
$\{ P < \lambda \} \subset \subset \Omega$ and we assume $(f/P)$ holomorphic outside a compact. 
In the case where $W g=g / Q$, that is $QW g=g$, 
where $Q$ is hypoelliptic, then $W$ corresponds to a
very regular distribution, that is $W \sim I$ modulo $C^{\infty}$.

\newtheorem*{res4}{Singularities}
\begin{res4}
 We assume all singularities for our model are on first surfaces to the symbol $f$ and can
 be reached through involutive movements, movements linked to movements orthogonal to involutive movements or algebraic
 approximations.
\end{res4}

\section{Invariant sets}
Consider $F(\gamma_{1} \rightarrow \gamma_{2})$ with $\gamma_{1} \bot \gamma_{2}$. Assume 
$\Phi : F(\gamma_{1})(\zeta) \rightarrow U_{1}$ and in the same manner 
$\Phi : F(\gamma_{2})(\zeta) \rightarrow \zeta \in U_{2}$. 
Assume $J \gamma_{1} = \gamma_{2}$. The problem is now under which conditions on $J$, 
do we have that there existence $\Phi$, such that $F(J \gamma)(\zeta) \rightarrow \zeta$ continuous. 
Sufficient for this is naturally that $F,J$ are analytic, why $\Phi$ is continuous. This is however not necessary. 
Assume $\gamma_{1}$ has a desingularization $U_{1}=\cup_{1}^{N} S_{j}$, where $S_{j}$ are connected. We 
can now define $J$ on the covering $J {}^{t} r_{j} \gamma_{1}={}^{t} r_{j} \gamma_{2}$ ($r_{j}$ is restriction). 
Note that when $\gamma_{2} \mid_{S_{1}} = \gamma_{2} \mid_{S_{2}}$, we do not necessarily have 
$\gamma_{2} \rightarrow \zeta$ uniquely. The condition ${}^{t} r_{j} J = J {}^{t} r_{j}$, 
means that $J$ is not 1-1. For every $\gamma_{2}$ partially hypoelliptic, naturally there is a $\gamma_{1}$ hyperbolic, 
such that $J \gamma_{1}=\gamma_{2}$

Monodromy $f(\zeta) \rightarrow \gamma$ is necessary in order to separate  $\gamma_{1}$ from 
$\gamma_{2}$. Assume $\gamma_{1} \rightarrow S_{j}$ and $\gamma_{2} \rightarrow \tilde{S}_{j}$ and 
$J : \gamma_{1} \rightarrow \gamma_{2}$, if we have $id : S_{j} \rightarrow \tilde{S}_{j}$ continuous, 
we can define $\gamma_{2}$ on $S_{j}$, that is on the covering to $\gamma_{1}$. More precisely, 
if $\widetilde{\Omega}$ is a covering defined by $\gamma_{1}$, we consider this as a domain for $\gamma_{2}$, 
that is $\{ \tilde{\gamma_{2}} \}=\cup (\tilde{S_{j}},\gamma_{2}(\tilde{S_{j}}))$ with analytic continuation and we write 
$\gamma_{2}(\widetilde{\Omega})=\tilde{\gamma_{2}}(\Omega)$. Concerning localization, starting from 
$F(\gamma) \rightarrow \gamma_{1} \rightarrow \gamma_{2} \rightarrow V$, if $V$ is accessible from 
$\gamma_{2}$, it is not necessarily accessible from $\gamma_{1}$, further accessible from 
$\gamma_{1} \rightarrow \gamma_{2}$ does not imply accessible from $\gamma_{2},\gamma_{1}$ or 
$F(\gamma)$. Further, if $F(\gamma) \rightarrow \zeta$ is not continuous, there are possibly $\tilde{\gamma}$, such that 
$F(\tilde{\gamma}) \rightarrow \zeta$ is continuous. 

Assume $F=\widehat{E}f$, where $E$ is a very regular parametrix to $\gamma_{2}$ and $F_{1}(\gamma_{1})=f$, then we must have that
$F(\gamma_{2})=f$, modulo $-\infty$ action. Thus, if $f \rightarrow f_{0}$, we have $\widehat{E}(f) \rightarrow f_{0}$
modulo $-\infty$. $E$ can be chosen with one-sided support, assuming $\gamma_{2}$ algebraic. When $E$
is constructed using Fredholm operators, $\gamma_{2} \sim \widehat{E}^{-1}(\varphi)$, as $\varphi \rightarrow \delta$. 
Note that $\{ \gamma_{2} \leq \lambda \} \subset \subset \Omega$, where $\lambda$ is a constant and $\Omega$ a domain
in $\mathbf{R}^{n}$.

Assume, for an analytic quotient, $(F_{2}/F_{1})(\gamma) \rightarrow 0$ over $\mid \gamma \mid=1$ with positive measure, 
we then have $F_{2} \bot F_{1}$ on $\mid \gamma \mid > 1$. More precisely 
$(F_{2} / F_{1})(r_{T}' \gamma)(\zeta)=(F_{2} / F_{1})(\gamma)(\zeta_{T})$ with $\mid r_{T}' \gamma \mid=1$. 
Define $E = \{ \zeta_{T} \quad \mid \gamma \mid=1 \quad (F_{2} / F_{1}) \rightarrow 0 \}$. 
Thus, if $(F_{2}/F_{1}) (\gamma) (\zeta) \rightarrow 0$ for large $\zeta_{T}$, that is $F_{1} \bot F_{2}$ 
with respect to $\zeta$ and $\mid \gamma \mid=1$, we then have $F_{1} \bot F_{2}$ with respect to $\gamma$.

Write $(F_{1} \bot F_{2})(\zeta)$, when the orthogonality is taken in $\zeta$, we then have 
$(F_{1} \bot F_{2})(\zeta)$ implies $(F_{1} \bot F_{2})(\gamma)$. In the same manner, if $F$ 
linear in $\gamma$, we have $\Delta(F)(\zeta)$ defines $\Delta(F)(\gamma)$.

Assume $f(\zeta_{T})=F(r_{T}' \gamma)(\zeta)$, where $r_{T}$ is assumed closed and locally 1-1, 
why $r_{T}'$ is locally surjektive. Define the continuation through $r_{T}'$ and 
$< \gamma_{2}, r_{T}' \gamma_{1}>=0$ or $< \gamma_{2}(\zeta_{T}),\gamma_{1}>=0$, note that the proposition 
that $r_{T}' \gamma_{1}=\gamma_{2}^{\bot}$ is equivalent with the proposition that $\gamma_{2}(\zeta_{T})$ 
is locally 1-1 at the boundary. Note that hyperbolicity assumes a Cartier boundary, while hypoellipticity 
assumes a bijective ramifier. 
We can define $r_{T}' (I) = (I)(\Omega_{T})$. 

Concerning accessibility, if for instance $\gamma_{1}=r_{T}' \gamma_{2}$, with $r_{T}' \rightarrow r_{T}$ 
locally 1-1 and closed (continuous), we have that $r_{T}'$ is surjektive $\sim J^{-1}$, that is 
$\forall \gamma_{1}$, we have existence of $\gamma_{2}$, such that $F(\gamma_{1} \rightarrow \gamma_{2})(\zeta)$ 
solves the lifting problem. We can assume that $J(\gamma_{1} \bot \gamma_{2})(\gamma)$ induces a continuous 
mapping $\tilde{J}(\gamma_{1} \bot \gamma_{2})(\zeta)$, that involves a complementary set to $\Omega$.

Assume $\gamma_{2}=J \gamma_{1}$, where we assume $\tilde{J} : N(\gamma_{1}) \rightarrow N(\gamma_{2})$. 
If we using Radon-Nikodym's theorem and let $F(e^{v} \gamma_{2})=F(J \gamma_{1})$, that is 
$F \circ e^{v} \sim F \circ J$, where $v \in L^{1}$. Then the condition $log J \in L^{1}$ implies 
algebraic singularities. $J$ is here defined as dependent on the movements, that is we write $v_{j}$ $j=1,2,3$.

Note that invariant sets are not necessarily preserved under iteration of symbols.
The sets $\Delta$ (lineality) and $\bot$ orthogonality can not be assumed independent for $(I_{1})^{\bot}$.
Consider for example $f=f_{1} + i f_{2}$ and
$\frac{f_{1} f_{2}}{f_{1}^{2} - f_{2}^{2}}=\frac{1}{\frac{f_{1}}{f_{2}} - \frac{f_{2}}{f_{1}}}=\frac{1}{v - \frac{1}{v}}$, 
implies $\frac{f_{1}}{f_{2}} \rightarrow 0$, when the quotient above $\rightarrow -0$ and 
$\frac{f_{2}}{f_{1}} \rightarrow 0 $, when the quotient above $\rightarrow +0$. Further let 
$e^{\phi}=f_{1}/f_{2}$, we then have $\mbox{ sinh } \phi=\frac{1}{2} \big[ e^{\phi} - e^{-\phi} \big]$, 
that has $\rightarrow \infty$ as $\phi \rightarrow \infty$ and $\rightarrow - \infty$ as 
$\phi \rightarrow -\infty$. Further, $1/\mbox{ sinh }\phi = \mbox{ csch }\phi$ that has 
$\rightarrow -0$, as $\phi \rightarrow -\infty$ and $\rightarrow +0$ as $\phi \rightarrow + \infty$.

For a hyperbolic symbol, we have that $f/Pr f$ is real (\cite{Garding87}). Let $f_{1}/f_{2} =  p_{m}^{1}/p_{m}^{2}$, where $p_{m}$ are highest order terms.
We then have that $p_{m}(tx)=t^{m}p_{m}(x)$, why if we choose $x$ such that $p_{m}^{1}(x) \neq 0$ 
and $p_{m}^{2}(x) \neq 0$, then we must have $f_{1}/f_{2}(tx)=const$, as $t \rightarrow \infty$,
why it is characteristic for hyperbolic symbols, that the real and imaginary parts are not orthogonal.

Assume $f^{N}=e^{\phi_{N}}$, with $\frac{d \phi_{N}}{d x} \rightarrow 0$ in the $\infty$, that is $d \phi_{N} \bot d x$
in $\infty$. 
The condition $\frac{1}{x} \phi_{N}(x) \rightarrow 0$ as $x \rightarrow 0$, implies that $\phi_{N} \rightarrow 0$
faster than $1/x$ goes to $\infty$, why the condition in $\infty$ implies a condition at the boundary. Assume $\phi_{N} + p_{N}=I$,
that is an algebraic complement. In this case we have $\phi_{N} \sim I - p_{N}$, that is invertible outside
constant surfaces to $p_{N}$, which gives possibility for two-sided limits. Constant surfaces for 
$\phi_{N}$ are constant surfaces to polynomials. An algebraic transversal implies oriented first 
surfaces (one-sided orthogonality). On the other side
if $d \phi_{N} + p_{N} d x=0$ then we have $\frac{d \phi_{N}}{d x}=-p_{N}(x)$ and when $p_{N}$ are 
reduced, we can assume $p_{N}(\frac{1}{x}) \sim \frac{1}{q_{N}(x)}$, for polynomials $q_{N}$.

Note the example $F(x,y)=\phi(x)\psi(y)$, where $\phi=\phi_{1} + i \phi_{2}$ and $\psi=\psi_{1} + i \psi_{2}$
We then have $F_{1} \bot F_{2}$ if $\frac{\psi_{2}}{\psi_{1}} + \frac{\phi_{2}}{\phi_{1}} / 1 - \frac{\phi_{2} \psi_{2}}{\phi_{1} \psi_{1}}$
that is if we have separately $\bot$ we have $\bot$.

Given analyticity we have $\tilde{\gamma} \in (I)(\Omega)=I(\tilde{\Omega})=\tilde{I}(\Omega)$, 
where we assume $\gamma \subset \tilde{\gamma}$ implies $\tilde{\Omega} \subset \Omega$. 
Further, $F \in \mathcal{D}_{\tilde{I}}' \subset \mathcal{D}_{I}'$, thus given that the continuation is allowed, 
we have that the restriction is well-defined. Assume invariance is defined by $\{ \gamma \quad F(\phi \gamma)=F(\gamma) \}$ 
for a transformation $\phi$. Define $\Omega_{(2)}=\{ \gamma \quad F(\phi \gamma^{2})=F(\gamma^{2}) \}$, 
it is then significant if $\phi \gamma^{2}=\gamma^{2}$ 
and $V_{2}=\{ \gamma \quad F^{2}(\phi \gamma)=F^{2}(\gamma) \}$. Consider $\Omega'=\{ \gamma \quad d F(\phi \gamma)=d F(\gamma) \}$ 
that is $\{ {}^{t} \phi d F(\gamma) = d F(\gamma) \}$ and $\Omega''=\{ d F^{2}(\phi \gamma)=d F^{2}(\gamma) \}$. 
Define $(J_{1})=\{ \gamma \quad {}^{t} \phi d F(\gamma)=d F(\gamma) \}$, for $F$ fixed. 
It is significant if ${}^{t} \phi$ is linear, when we study $\Omega=\Omega_{1} \cap \Omega_{2}$. 
$(J_{1}) \subset (J_{2})$ implies $\Omega_{2} \subset \Omega_{1}$ according to algebraic geometry 
and $(J_{2})' \subset (J_{1})'$ according to functional analyse. 
We can if $F^{2}$ is linear over $d \gamma$, 
conclude $\phi d \gamma - d \gamma \in \mbox{ ker }F^{2}$. If $F(0)=0$, we have $\mbox{ ker }F \subset \mbox{ ker }F^{2}$. 

Assume that $\Delta$ defines a geometric ideal $(I)$, that is $g \in (I)$ implies $\tau f - f = g$, 
for some $f$. 
Assume that $\Omega_{j}$ gives invariant sets, then we have that $\Omega_{1} \cap \Omega_{2} \rightarrow (I_{1}) + (I_{2})$. 
Sufficient for a disjoint decomposition is that $(I_{j})$ are given by positive functions. 
Consider $f \in (I_{1})$ and $f=f_{+} - f_{-}$. We have that $\{ 0 \} \subset \Omega_{1} \cap \Omega_{2}$ 
implies $(I_{1}) + (I_{2}) \subset (J)$, where $(J)$ is the ideal corresponding to a disjoint decomposition.

Assume $f_{2}/f_{1}=\varphi$ and $M(\phi) \sim \varphi$ (arithmetic mean), then the proposition that $\phi$ is harmonic corresponds to, 
$\frac{d}{d z} \frac{f_{2}}{f_{1}}$  real and analytic. When $\frac{d}{d z} \frac{f_{2}}{f_{1}}$ real, 
we have $\frac{f_{2}}{f_{1}} \frac{d}{d z} \log \frac{f_{1}}{f_{2}}=\frac{d}{d z} \log \log \frac{f_{1}}{f_{2}}$.
Thus, if $\log \varphi \in L^{1}$ and $\frac{d}{d z} \log \log \frac{f_{1}}{f_{2}}$ 
is real and analytic, then we have existence of $\phi$ harmonic, such that $M(\phi) \sim \varphi$ 
is constant. Define $F(e^{\varphi})=\widehat{F}(\varphi)$ and $(\mathcal{F} S)(\varphi)=S(e^{\varphi})= S \exp \varphi$ 
and $(\mathcal{F}^{-1}S)(\varphi)=S \log \varphi$. We then have 
$\widehat{(\log \frac{f_{2}}{f_{1}})}=\mathcal{F} I \log \frac{f_{2}}{f_{1}}=\mathcal{F} \mathcal{F}^{-1} I \frac{f_{2}}{f_{1}} \sim \frac{f_{2}}{f_{1}}$
Thus the condition $\log \frac{f_{2}}{f_{1}} \in L^{1}$ using the Lelong transform, 
can be interpreted so that $f_{1} \bot f_{2}$ (one-sided).

Note that in the model in this article, an algebraic continuation of an hyperbolic operator 
is not hyperbolic (and vice versa). Note that if the ideals are given by distance functions $d_{1},d_{2}^{\bot}$,
such that $d_{2}^{\bot} / d_{1} \rightarrow 0$ in $\infty$ and $d_{2}^{\bot}/d_{1} \leq \mid \tilde{\gamma} \mid$,
where $\tilde{\gamma}=J \gamma$, which does not imply $d_{1} \leq \mid \gamma \mid$.
Existence of $N$, such that $P$ is hyperbolic in the direction $N$, implies that $P$ is not hypoelliptic. 
Since the principal part $p_{m}$ to a 
phe operator, is independent of some variables, we have that there exists $N$ where $P$ is hyperbolic 
(in dependent variables). Thus the restriction of $(I_{Phe})$, can in this way be a subset of $(I_{Hyp})$. 
More precisely, given $(I)=(I_{Hyp}) \oplus (I_{Phe})$, when $(I_{Hyp})$ is seen as an extension of 
$(I_{Phe})$, there is a corresponding restriction to dependent variables for $p_{m}$, as the domain 
for $(I_{Hyp})$. Consider $\frac{1}{t^{k}} f(t x) \rightarrow Pr f(x)$, as $t \rightarrow 0$, 
where $k$ is the order  of zero. Particularly $t^{-N} f(t \eta) = Pr f(0)$ $\forall N$ and 
$t \eta \in \Delta$, that is an ``infinite zero''. For a reduced operator, we have thus 
$t^{N}f(x/t) \nrightarrow Pr f(0)$, for $\forall N$ and $t \rightarrow \infty$. 
Reduced operators are of real type (type 0) and modulo $C^{\infty}$ we can  always assume real type.

Assume that $S_{1}^{\bot}=\{ $ normals to $S_{1} \}$ $=N(f_{1}^{\bot} + f_{2})$, where $S_{j}$ are first surfaces. 
The condition $S_{2} \cap S_{1}^{\bot} \neq \emptyset$ implies $S_{2}^{\bot} \cap S_{1} \neq \emptyset$. 
Further $(f_{2} + f_{1}^{\bot})^{\bot}=f_{2}^{\bot} + f_{1}^{\bot \bot}$. The mapping $S_{2}^{\bot} \cap S_{1} \rightarrow S_{2} \cap S_{1}^{\bot}$ 
is a mapping between layers. 

Consider the first surfaces to the iterated symbols and $M$ a retraction neighbourhood of $S_{j}$.
If we assume that $S_{j}$ stratum, we assume an embedding $S_{j} \hookrightarrow M$, 
that is given a stratification $\{ S_{j} \}$, we assume 
$S_{1} \hookrightarrow M \rightarrow_{projection} S_{2}$, which gives a mapping between $S_{1}$ and $S_{2}$.
Assume $S_{1}^{2} \supset S_{1}$ and assume the same conditions for $S_{2}$, consider 
$S_{1}^{2} \rightarrow S_{1} \hookrightarrow M \rightarrow S_{2}^{2} \rightarrow S_{2}$, 
and assume the compatibility condition, $S_{2}^{2} \cap S_{2} = \emptyset$ or $S_{2}^{2} \supset S_{2}$.
The condition that $f_{2}/f_{1}$ is polynomial in the infinity, is necessary to come to the conclusion
that the intersection $\bot$ is discrete (for instance \cite{Julia23}). The condition is also necessary to conclude that 
the condition $\bot$ is independent of $\zeta$ when $\mid \zeta \mid \geq R$, $R$ large. When $\Delta_{(2)} \subset \Delta_{(1)}$ 
(index refers to iteration) we have for the first surfaces that intersect $\Delta$, that 
$\tilde{S}^{1} \subset \tilde{S}^{2}$.
In the two systems we assume the lineality is the same,
it is characteristic for $(I_{1})$ that $S_{1} \sim S^{2}_{1}$. According to the compatibility conditions, 
assuming $S \rightarrow \tilde{S}$
we can assume $\tilde{S}^{1} \subset \tilde{S}^{2}$ or $\tilde{S}^{1} \cap \tilde{S}^{2} = \emptyset$

Assume $\gamma_{1} \rightarrow \gamma_{2}$ with the corresponding first surfaces $S_{j}$ and 
$\tilde{S}_{j}$, where we assume $\gamma_{1} \rightarrow \gamma_{2}$ continuous and $\gamma_{j}$ analytic, 
that is we do not necessarily have existence of $S_{j} \rightarrow \tilde{S}_{j}$ continuous. 
For the problem of continuing a desingularization to $\gamma_{2}$, it is sufficient to assume $\gamma_{1} \rightarrow \gamma_{2}$
locally 1-1. Alternatively, if we define $S_{j}$ through a distance function $d_{1,j}$ and in the same
manner for $\tilde{S}_{j}$ and $d_{2,j}$, we can define $d_{1,j}d_{2,j}$ and a ramifier $r_{T}'$ 
from $S_{j}$ to $\tilde{S_{j}}$. Injectivity for $r_{T}'$ 
is regarded relative $d_{1}d_{2}=0$. We assume Schwartz type topology according to (\cite{Martineau})), inclusion
of ideals $I_{1,j} \subset I_{1,j+1}$ is defined through $d_{1,j+1}/d_{1,j} \rightarrow 0$.
This implies surjectivity for an approximating sequence.  
Thus, the condition $\gamma_{1} \rightarrow \gamma_{2}$
continuous, linear and locally 1-1, is sufficient for existence of path $S_{j} \rightarrow \tilde{S}_{j}$.
 
\newtheorem*{res5}{Continuation by continuity}
\begin{res5}
 A sufficient condition on $J : \gamma_{1} \rightarrow \gamma_{2}$ continuous, to induce a continuation
 of the correspondent first surfaces, is that $J$ is continuous, linear and locally 1-1.
\end{res5}
Note that this property is not necessary for a global mixed model, that is $f(\zeta) \rightarrow \gamma_{1} \rightarrow \gamma_{2} \rightarrow V$,
for a given geometric set $V$

\section{The two mirror model}
We discuss the mixed problem in a two base (two mirror) model, that is we consider $f(\zeta)=F(\gamma)(\zeta)$,
where $\gamma \in (I) \bigoplus (I)^{\bot}$. Assume $J_{\delta} \gamma_{1} \rightarrow \gamma_{\delta} \in (I_{1})^{\bot}$,
for a parameter $\delta \rightarrow 0$,
then given $F(V \gamma_{\delta}) - F(V^{*} \gamma_{\delta}) \sim 0$, where $V$ is a given movement,
the limit $V \rightarrow I$ exists as two sided. Further, When $F \sim GH$ and $F(\gamma_{1},\gamma_{2}^{-1}) \sim 
\int G(\gamma_{1},\mu) H(\mu,\gamma_{2}^{-1}) d \mu$, when $\mu \rightarrow \delta_{0}$, then the two 
mirror model goes over into a one mirror model and we have $F \sim (G \otimes I)(I \otimes H)$. Assume 
the boundary $\Gamma$ can be given by one single function $\eta=y(x)/x$ (order 0). Denote $\eta_{1}^{*} \sim V^{*} \eta_{1}$ 
reflection through $V=I$, then we have in the two mirror model, $(\eta_{1}^{*})^{*} \sim \eta_{2}^{-1}$.
Assume $< U \eta_{1},\eta_{2}>=< \eta_{1}, V \eta_{2}>=0$ and $\mu \eta_{1} = U^{*} \eta_{1} \sim V^{-1} \eta_{2}$
and we assume $\eta_{1}=\eta_{2}$, over $\mu$ and $\Delta$. Then $\mu \eta_{1} \rightarrow \eta_{1}$
and $\mu^{-1} \eta_{2} \rightarrow \eta_{2}$, that is ${}^{t} \mu = \mu$ is involutive. When 
$\mu \rightarrow \delta_{0}$ we have $\eta_{1}(0)=\eta_{2}(0)$. Otherwise, we assume $\mbox{ supp }\mu=\Delta$. 
When $F$ is symmetric over the path, we have that $V$ is involutive relative $F$. 
In the system $(\eta_{1},\eta_{2})$, when we consider $\Delta \rightarrow \mu$,
then $ \eta_{1}=\eta_{2}$ can be regarded as an abstract light cone, that is $\eta_{1} \times \eta_{2} \sim \mu$.

Concerning the
boundary condition $F(\eta_{2}) \mid_{t=0}=F(\eta_{1})$ and $\frac{d F}{d t}(\eta_{2}) \mid_{t=0}= \frac{d F}{d t}(\eta_{1})$. 
In particular, when $\frac{d F}{d t}= \frac{d F}{d x} \frac{d x}{d t} + \frac{d F}{d y} \frac{d y}{d t}= -Y_{2} X_{1} + Y_{1} X_{2}=0$,
we have $\frac{Y_{2}}{X_{2}}=\frac{Y_{1}}{X_{1}}$ according to the Lie's point transform. 
For the lifting operator we use an involutive condition,
that is if $\sharp{F}$ is the continuation of $F$ according to $\gamma_{1} \rightarrow \gamma_{2}$,
we assume $\frac{d F}{d x} \frac{d \sharp{F}}{d y} - \frac{d F}{d y} \frac{d \sharp{F}}{d x}=0$
in a neighbourhood of $\Delta$, which implies Lie's condition as above.

The lineality for the composite kernel can be represented as $\big[ G,H \big] (\xi + i t \eta)=\big[ G,H \big]$
and $\int G(x,y)H(y,z + i t \eta) d y=\int G(x + i t \eta,y)H(y,z) d y$. Thus, if $\big[ G,H \big]=\big[ H,G \big]$,
then $\Delta(\big[ G,H \big])=\Delta(G) \cap \Delta({}^{t}H)$. If $G$ has $\Delta(G)=\{ 0 \}$, the same holds for 
the composition, only assuming symmetry.
Further that if $G=G_{1}G_{2}$ with $G_{1}$ hypoelliptic and $N(G_{1}) \neq N(G_{2})$, we do not necessarily have that $G$ is hypoelliptic, 
since $\frac{G'}{G}=\frac{G_{1} '}{G_{1}} + \frac{G_{2} '}{G_{2}}$, where only one term in 
the right hand side is assumed $\rightarrow 0$.

\newtheorem*{res6}{Two mirror movements}
\begin{res6}
 The two mirror model refers to an involutive movement, using two reflection points. A generalization 
 to the use of two bases $\gamma_{j}$, requires localizers $F_{j}$ with one-sided support $j=1,2$. A necessary
 condition for the two mirror model, to give a normal model, is that one of the limits is independent of orientation.
\end{res6}

Given a reflection in the plane with reflection points on the axles, assume that
the distances from $B$ to $A$ is given by $d_{1},d_{2},d_{3}$. We then have 
$d_{3} \rightarrow 0$, $z \rightarrow A$ and $\rightarrow \infty$
as $z \rightarrow A'$. Further $d_{2} \rightarrow 0$ as $z \rightarrow A'$, $\rightarrow \infty$ 
as $z \rightarrow B'$.
Finally $d_{1} \rightarrow \infty$ as $z \rightarrow B$ and $\rightarrow 0$ as $z \rightarrow B'$, 
that is $\frac{1}{d_{1}} \rightarrow 0$
as $z \rightarrow B$ and $\rightarrow \infty$ as $z \rightarrow B'$. 
Note that $\frac{1}{d}$ does not necessarily define a distance.

\subsection{Orientation of limits}

Given a simply connected domain $\Omega$ in the plane, then 
we have that a simple Jordan curve $\Gamma$ divides $\Omega$ into two connected components. 
Assume $\Gamma$ is a simple curve that is transported in a normal model along the transversal. 
If $V_{1},V_{2}$ are two sets intersecting the transversal, considered as ``antipodal'', 
we do not necessarily have existence of a unique plane between $V_{1},V_{2}$, that intersects
the transversal. However according to Nishino (\cite{Nishino68}), if $V_{1},V_{2}$ are first 
surfaces with constant values $c_{1},c_{2}$ and if $c_{1} < c_{2}$, we can determine the intermediate set 
as a first surface to a scalar $c$, $c_{1} < c < c_{2}$. 
Assume $d_{j}$ distances to respective set and 
consider $< d_{j}, \mid f \mid >$, then if we assume $d_{1} / d_{2} \rightarrow 0$ 
in $\infty$, we have a continuous injection between respective spaces, weighted with 
the distances. The relation defines a topological inclusion between respective ideals.
More precisely, under the condition $d_{1}/d_{2} \rightarrow 0$ in $\infty$, we can form 
$L_{d_{2}},L_{d_{1}}$, 
that is $L_{d_{2}} \subset L_{d_{1}}$. Thus, $J \gamma_{1} \in L_{d_{2}}$ $\Rightarrow \gamma_{1} \in L_{d_{1}}$. 
Further $\int d_{1} \mid J^{-1} f \mid \leq \int d_{2} \mid f \mid$. 

Assume $\eta=y/x \in \dot{B}$, we then have existence of $F \in \mathcal{D}_{L^{1}}'$ 
such that $F \sim P(D) \tilde{f}_{0}$, where $\tilde{f}_{0}$ very regular. When we have 
$\eta \in B$, we assume existence of $\eta_{1} \in \dot{B}$ (\cite{Dahn13}), such that $\eta \sim_{m} \eta_{1}$. 
Note that when $\eta_{2}$ corresponds to $\gamma^{2}$, if we have $\eta \in \dot{B}$ or $1/\eta \in \dot{B}$,
we have $\eta_{2} \in \dot{B}$. Existence of limits in $\mathcal{D}_{L^{1}}'$ can be seen 
as a topological type (A) condition (\cite{Nishino75}). Note that when we do not have existence of limits $F(\frac{1}{\eta_{2}})$ 
in $\mathcal{D}_{L^{1}}'$, we have existence of limits $F(\eta_{2}^{\bot})$ (annihilators).

Consider $F(\gamma_{1} \rightarrow \gamma_{2}) \rightarrow U_{2}$ and $F(\gamma_{1} \rightarrow \gamma_{3}) \rightarrow U_{3}$
given $U_{2}^{\bot} \supset \Delta (\gamma_{j})$, $j=1,2$. When $r_{T} \zeta$ is locally 1-1 and 
closed, we have that $r_{T}' \gamma$ is surjective. Thus, for every $\gamma_{2}$ we have existence of
$\gamma_{1}$ such that $J \gamma_{1}=\gamma_{2}$, thus starting from $\Delta$ we can always find 
a homogeneous symbol with $\Delta$. Further, when $F(\gamma_{1} \rightarrow \gamma_{2}) \rightarrow U$
and $F(\gamma_{1} \rightarrow \gamma_{3}) \rightarrow U$, then $\gamma_{2}=r_{T}' \gamma_{3}$, 
that is we have existence of $\mu$, a continuous path between $\gamma_{2}$ and $\gamma_{3}$. Thus the fact that we have 
existence of an approximation of $U$, does not exclude existence of a longer path with the same limit.
In applications, the paths may have different order of zero's, thus different quality properties.
Starting from the moment problem and $J \gamma_{1}=\gamma_{2}$, given $\gamma_{2}$ a polynomial, we 
have that $\liminf \gamma_{2} \leq \gamma_{2} \leq \limsup \gamma_{2}$, for the restriction to lines. For instance when $f$ is 
of type (A), as long as $F(\gamma_{1} \rightarrow \gamma_{2})$ preserves this inequality and with finite limits,
we can solve the moment problem.  When the limits coincide, the solution is unique.

\newtheorem*{res7}{Twosided limits}
\begin{res7}
 Given $\gamma_{1}$ from $A$ to $A'$ (reflection point), $\gamma_{2}$ from $B$ to $B'$ and $A' \sim B'$.
 We refer to this as a two-sided limit. When $\gamma_{2}$ also gives a path from $B'$ to $B$, the composition
 gives a path from $A$ to $B$.
\end{res7}

We can divide two mirror model into $RF(z)=F(\overline{z})$ and $SF(z)=F(i z)$, that is we assume the reflection
points on the real axes or the diagonal. 
When we do not have symmetry according to $F(\overline{z})=\overline{F}(z)$, 
we must consider $f(z,\overline{z})$. In the same 
manner if $F(i \gamma)=i F(\gamma)$ or rather $(F + i F^{*})^{*}=-i (F + i F^{*})$ (pure), we can 
refer to this as a transmission property, otherwise we must consider $F(z,iz)$ (or $F(w,w^{*})$).

Assume in the two-mirror model, that $L$ is the segment between $\gamma_{1}$ and $\gamma_{2}$ and 
consider $F(\gamma_{1} \rightarrow \gamma_{2})=\int G(\gamma_{1},L) H(L,\gamma_{2}) d L =\big[ G,H \big](\gamma_{1},\gamma_{2})$. 
Given $F$ invariant for change of orientation, this corresponds to one-sidedness (with respect to two 
mirrors). If $L \bot \gamma_{1}$ and $L \bot \gamma_{2}$, that is defined for instance by the distance 
functions, we then have existence of $L \sim \gamma_{1} \times \gamma_{2}$,
where $L$ is assumed to have points in common with $\gamma_{1},\gamma_{2}$. 

\subsection{Wellposedness}

Every movement in $(I_{Phe})^{\bot}$ can be represented as reflection (\cite{Riesz43}). This means that we have that 
$\gamma_{1} \in (I_{Phe})^{\bot}$ and $F(\gamma_{1})=const$, can be reached through $U \gamma$ with 
$\gamma \in (I_{Phe})^{\bot}$ and given by reflection relative an axes $R$ (not unique), that is $U=U_{R}$. 
Assume compatibility conditions, $F(U^{\bot} \gamma) \neq const$, for the 
approximation and assume that there are points in common, but $U^{\bot} \gamma$ is not necessarily in $(I_{Phe})^{\bot}$. Given that $\gamma \in \Sigma(S)$, a first surface, 
we have that the  reflection axes is a part of $\Sigma(S)$. Thus, if $R^{\bot}$ is the 
reflection axes to $U^{\bot}$ as (an euclidean) $\bot$ to $R$, it can be used as a regular approximation.

F. and M. Riesz theorem: assume $f(\frac{1}{z})$ analytic and bounded with 
$\mid f(\frac{1}{z}) \mid <M$ for $\mid \frac{1}{z} \mid < 1$. 
Let $E=\{ f(e^{i \theta})=\lim_{r \rightarrow 1} f(r e^{i \theta})=0 \}$. 
Assume the measure for $E$ is $>0$, on $\mid z \mid=1$, we then have  
$f \equiv 0$ for $ \mid \frac{1}{z} \mid < 1$ (\cite{Riesz56}). 
Assume $\Sigma= \{ \gamma \quad F(\gamma)=const \}$ equipped with a norm $\Sigma_{\rho}=\{ \gamma \quad \rho(\gamma)=1 \}$.
Application on my model, gives that if the segment $\mu$ exists
between $\gamma_{1}$ and $\gamma_{2}$ in the boundary $\Gamma$, forming a set of positive measure,
and if two solutions $F_{1} \equiv F_{2}$ on $\Sigma_{\rho}$, we have that $F_{1} \equiv F_{2}$ on $B_{\rho}=\{ \gamma \quad \rho(\gamma) \leq 1 \}$. 
In this application, the boundary is continued with a transversal between the first surfaces. The continuation
principle (\cite{Martineau}) gives a representation of $F$ as a distribution of the continuation.

Assume $J_{1}=J_{x} \otimes J_{y}$ and assume $< (J - J_{1}) \gamma, T \gamma >=0$, $\forall \gamma$, implies $T$ is the id-mapping,
we then have that $J=J_{1}$. Note for $\gamma$ analytic, we have $\gamma=0$ iff $\gamma \mid_{L} = 0$
for every line $L$, that is we can consider a pluri complex formulation. 
Note in this case if $\{ J \gamma_{1} < \lambda \}$ is semi algebraic locally, 
we have that for instance $\{ J_{x} \gamma_{1} < \lambda \}$ is semi algebraic.

\subsection{Localisation}
Assume $\eta(x)=y(x)/x$ exact, then the ideal given by $\mbox{ ker }y$, can be represented with a
global pseudo base. The same proposition holds if $y$ is reduced. 
An ideal that is defined by $\tau_{z} \phi/\phi \sim_{m} 0$
with $\tau_{z}$ compact, can be given a global pseudo base and when $\eta>0$ over an ideal, then 
$\eta$ is exact and the ideal has a global pseudo base. 

\newtheorem*{res8}{Localisation problem}
\begin{res8}
 Given $F(\gamma_{0})(\Omega)$, determine $\gamma_{2}$ analytic 
on $V$, such that $F(\gamma_{0}) \rightarrow \gamma_{0} \rightarrow \gamma_{2}$ continuous, further $\gamma_{2} \rightarrow V$ 
continuous and $\gamma_{2}^{N} \in (I_{He})$, that is representation using reduced measures.
\end{res8}
Assume $\{ \Omega_{j} \}$ a covering to $\gamma_{0}$ and $\Omega$. Given a point $\zeta_{0} \in V$
and $V_{k}$ a component, if $\gamma_{0} \in (I_{He})^{\bot}$, then we can give $\gamma_{2}$ as a 
continuation on components $\Omega_{j}$ to $V_{k}$ and $< \gamma_{0},\gamma_{2}^{N} >=0$. When we have 
$d \mu_{2} \in \mathcal{E}^{'0}$, $\gamma_{0}$ can be extended with zero over $V_{k}$. 

Assume $F \gamma_{2}=1$ (invertible), where $\gamma_{2}$ is a hypoelliptic symbol, 
and $F$ has representation with trivial kernel. 
Assume $\mbox{ ker }F e^{v}=\{ 0 \}$, where $F(e^{v} \phi)=0$ 
implies $\phi=0$ or $e^{v} \phi \in \mbox{ ker }F$ for $v \in L^{1}$. We assume 
$F(J \gamma_{1})=F(e^{v_{1}} \gamma_{2})=F(e^{v_{1} + v_{2}} \gamma_{2}^{2})$.
Thus, if $\gamma_{2}^{2}$ is hypoelliptic, we have $e^{v_{1} + v_{2}} \rightarrow 0$ in $\infty$,
that is $Fe^{v_{1} + v_{2}} \rightarrow 0$ and $F \nrightarrow 0$ in $\infty$. Then 
$\log \mid F \mid + v \leq \log \mid F \mid$, we have that the type for $v \leq 0$, that is of type $-\infty$. 
Note that if $JU=VJ$, then $U,V$ do not behave algebraically similar, that is $U \gamma_{1}^{2} = \gamma_{1}^{2}$ 
iff $U \gamma_{1}=\gamma_{1}$, but when $J \gamma_{1}=\gamma_{2}$, we do not have $J \gamma_{1}^{2}=(J \gamma_{1})^{2}$, 
that is $J$ is not algebraic.

Note that if $f(x)=\frac{1}{x}g(\frac{1}{x})$ and $f(\frac{1}{x})=\frac{1}{h(x)}$, we have that 
$x=\frac{f(\frac{1}{x})}{g(x)}$, that is $f(\frac{1}{x})h(x)=1$ and $f(\frac{1}{x}) \frac{1}{g(x)}=x$.

\subsection{Algebraic approximations}

The mapping $F(\gamma) \rightarrow \gamma \rightarrow \gamma_{2} \rightarrow \zeta$, does not necessarily preserve the 
order of $0$ for the localisation. For the localisation, existence of order $0$ regular approximations is 
sufficient. Note that concerning $\tau \rightarrow V_{e}$ using a contact transform, we have that the order of $0$ 
is preserved, but not the shape of obstacles. We will assume the boundary locally of order $0$ 
for some movement, that is it can be defined by a orthogonal movement locally. Assume $F(U^{\bot} \gamma) \rightarrow c$ 
regular and $F(U \gamma) \equiv c$, 
given that $\gamma$ is polynomial, we have through the lifting principle, existence of $({}^{t} U^{\bot} F)$ with an analytic representation. 
It is sufficient, that ${}^{t} U^{\bot} F \in \mathcal{D}_{L^{1}}'$ with $({}^{t} U^{\bot} F)(\gamma)(\zeta)$ 
analytic, that is $F \sim F(\zeta,\vartheta)$ kernel in $\mathcal{D}_{L^{1}}'$.

Assume $\{ \zeta \quad F(\gamma_{1})(\zeta)=c_{1} \}=S_{1}$ 
and $(I_{c_{1}})=\{ \gamma_{1} \quad F(\gamma_{1})=c \}$. Assume $\frac{d F}{d x_{1}}=-Y_{1}=0$ 
in $\zeta$ with $S_{1} \cap \{ Y_{1}=0 \}$ algebraic and the same condition for $X_{1}$. Assume $S_{1} \rightarrow \tilde{S}_{1}$ a continuation 
by continuity. Assume $J : \gamma_{1} \rightarrow \gamma_{2}=\gamma_{1}^{\bot}$, where $\gamma_{1}^{\bot}$ 
is closed, given that we have existence of $F^{-1}$ continuous, we have $F \rightarrow \gamma_{1} \rightarrow \gamma_{1}^{\bot}$.
As $\{ F-c \} \rightarrow 0$, we have $\gamma_{1}^{\bot} \rightarrow {}^{o}(I_{c})$ (annihilator). When $\gamma_{1}^{\bot}$ 
partially hypoelliptic, we have $(\gamma_{1}^{\bot})^{N}$ locally 1-1 (downward bounded) which implies 
$\zeta \rightarrow \zeta_{0}$, that is the limit exists. 

For instance assume $\gamma_{1} \in (I_{c})(S_{1})$ and $\gamma_{2}$ such that $\mbox{ supp } \gamma_{2} \cap S_{1} = \emptyset$ 
and $F(\gamma_{2})=const$. If $(I_{1}) \subset (I_{He})^{\bot}$ 
where $(I_{1})=\{ d_{1}=0 \}$.
Note that $d_{2}^{\bot} / d_{1} \rightarrow 0$ defines a continuation of $\gamma_{1}$, that is 
$\tilde{(I_{c}})=\{ \gamma \quad F(\gamma)=c \quad \gamma \bot \gamma_{2} \}$ and the ``moment problem'' 
gives that $(I_{1})$ has a continuation to $(I_{1})^{\bot}$, assuming
algebraic singularities. 

Assume $J_{\delta} : \gamma_{1} \rightarrow \gamma_{2}$, for a parameter $\delta \rightarrow 0$, 
such that $({}^{t} J_{\delta} F)(\gamma_{1})=c$ over an involutive set, that is $J_{\delta} \gamma_{1} \in F^{-1}(c) \sim (I_{c})$. 
Given type (A), it is sufficient to assume $J_{\delta} \gamma_{1}$ a polynomial. If $\gamma_{1}^{\bot} \in \mathcal{D}_{L^{1}}'$ 
which implies $\gamma_{1}^{\bot} \sim P \mu$ and $P \mu (\gamma_{1})=\mu ({}^{t} P \gamma_{1}) \sim$ 
algebraic continuation close to the boundary.

Assume the continuation is given by a movement, 
such that $(U^{\bot} \gamma)$ and $(U \gamma)^{\bot}$ have points in common.
For instance $U^{\bot} \gamma \rightarrow \{ F(\gamma)=c \}=S_{1}$ regularly and $F(U \gamma)=c$ 
and $\tilde{U} \zeta \subset S_{1}$. Let $V^{\bot} \gamma=J_{\delta} U^{\bot} \gamma_{1}$, where we assume
$\{ V^{\bot} \gamma \} \cap \{ F(\gamma_{2})=c \}=\tilde{S}_{1}$ algebraic, this can be seen as a 
compatibility condition (when $v \in L^{1}$).

Consider for example $F$ as an analytic function of $\gamma$, where the boundary $\Gamma$ is defined as $\bot$ 
movements, that is if $V \gamma$ is translation, and $V=JU$, we can assume
$U^{\bot} \gamma$ are parallel planar movements. Given a normal model, if the measure for 
$V \gamma$ is positive and $F(V \gamma)=f(\zeta)$, we have wellposedness according to a previous argument. 
Note $F(V^{\bot} \gamma)= const$, that is $U \gamma$ denotes a $\sim$ planar movement. 
If $R_{1}$ is the reflection axes and $R_{2}$ the corresponding axes in euclidean geometry. 
Assume $\tilde{R_{2}}$ the corresponding points to $\tilde{V}$, that is $V \gamma = \gamma$ implies
$\zeta \in \tilde{R_{2}}$, then if $V-I$ has an analytic representation, this is a line.

\section{Discussion on change of base}

Assume $f \in (I)$ and $g \in (I^{\bot})$, we define $H(g)=H(f)=0$ ($H=\tau-1$) on $\Delta$, that is
$H(f)=A_{f}(H)$ and $H(g)=B_{g}(H)$ which implies $A_{f}(H)=B_{g}({}^{t} e^{\varphi} H)$ as 
``inverse functionals''. Thus $H(g)=H(e^{\varphi} f)$, where $\varphi \in L^{1}$ in our model, which 
corresponds to equivalent zero's. Consider $F(f+g) \sim F_{2}((f - e^{\varphi} g) + i (g + i e^{\varphi} f))$ 
assuming $F_{1} \bot F_{2}$. Assume $R(f)=e^{\varphi} g$ and $R^{*}(g)=-e^{\varphi} f$ where $J$ is reflection through $\Delta$ 
and $F(J f)=F(e^{\varphi} g)$, that is $F \sim F_{2}((1-R)(f) + i (1-R^{*}(g)))$. 

Starting with $V \subset \Omega$, and $\Omega \backslash V$ analytic, we can  form $(I)(V) \otimes (I)(\Omega \backslash V)$
Denote by $(I^{\bot})=\{ g \quad fg = 0 \}$, an annihilator ideal. The corresponding geometric set
is $V \cup \Omega \backslash V$. 
If we assume $f,g$ positive, we have
$N((1-f)(1-g))=N(1 - fg)$, that is $(fg)^{\bot}=f^{\bot} g^{\bot}$.
Assume instead $(I) \oplus (I^{\bot})$ with base elements $f,g$ and $F_{1}={}^{t}e^{\varphi}F_{2}$, where ${}^{t} \varphi \rightarrow \infty$.
We then have $F(\alpha f + \beta g)=F_{2}((1 + i e^{\varphi})(\alpha f + \beta g))=F_{2}(\alpha f - e^{\varphi}  \beta g + i (\beta g + e^{\varphi} \alpha f))$
and $= F_{2}( \alpha f - e^{\varphi} \beta g) + i F_{2}(\beta g + e^{\varphi} \alpha f)$.
Further, if $\alpha=\beta=1$, $\frac{g + e^{\varphi} f}{f - e^{\varphi} g}=
\frac{\frac{g}{f} + e^{\varphi}}{1 - e^{\varphi} \frac{g}{f}}$, why if $\frac{g}{f} \rightarrow 0$ in $\infty$
we have that the quotient above $\rightarrow \infty$. In the same manner, if $f/g \rightarrow 0$  in $\infty$
faster than $e^{\varphi} \rightarrow \infty$, the quotient $\rightarrow -0$ in $\infty$. 
Thus, we have that $F_{1} \bot F_{2}$ one-sided, given $f \bot g$ one-sided.

Assume $(I)^{\bot}=\{ g \quad <f,g>=0 \}$ on a set of positive measure, then we have that $\mbox{ supp }g=W$
is the set where $f$ is not a polynomial. That is, $\int_{W} fg d x=0$ implies $fg \equiv 0$ on $W$ or that the measure for $W$ 
is zero and if $W$ contains a connected set where $fg \equiv 0$, we have that this set has measure zero.
If $fg \equiv 0$ on a line $L$ $\subset \Delta$ , this implies $g$ a zero divisor. Assume $\tau'=\tau-1$ 
algebraic and $\tau'(fg)=\tau' f \tau' g$ with $\tau f = f$ and $\tau' (fg) \equiv 0$, independent of 
$g$. Assume now $< \tau' f,g > \equiv 0$ and $< f,\tau' g> \neq 0$. Define 
$\Sigma=\{ < \tau' f,g >=< f,\tau' g> \}$. 
If $g$ is reduced, we thus have $\Sigma = \{ 0 \}$. 
Assume $(J_{f})=\{ g \quad <\tau' f,g >=< f,\tau' g> \quad \mbox{ on } L \}$ where
$L$ has positive measure. We then have $\tau' - {}^{t} \tau'$ is not algebraic! 
We have that $L \subset N(J_{f})$ and $I N (J_{f}) \sim rad (J_{f})$.

Assume $\gamma_{1} \rightarrow \gamma_{2}$ defines a broken ray according to $\gamma_{2}^{* -1} \simeq \gamma_{1}^{*}$
thus $\gamma_{1},\gamma_{2}$ must have a point in common. It is sufficient for this to assume $\gamma_{2}$ 
has two sided limits.
Assume $(\gamma_{2}^{-1})^{*} \simeq (\gamma_{2}^{*})^{-1}$ which means $\gamma_{1}^{*} . \gamma_{2}^{*} \sim 1$
(Legendre). We denote this $\gamma_{1}^{*} \mathcal{L} \gamma_{2}^{*}$, $\gamma_{1} \bot \gamma_{2}$,
$\gamma_{1} \mathcal{L} \gamma_{1}^{*}$, $\gamma_{2} \mathcal{L} \gamma_{2}^{*}$, $\gamma_{1}^{*} \mathcal{L} \gamma_{2}^{*}$,
that is we can define $\bot \rightarrow \mathcal{L}$ continuous. Note that $x . x^{*} + y . y^{*} \sim 1$ 
assumes $x \bot y$, that is when $\gamma_{1}^{*} \mathcal{L} \gamma_{2}^{*}$, then we must have $\gamma_{1} \bot \gamma_{2}$. When reflexivity is given by $\gamma_{1}^{**}=(x_{1}^{**},y_{1}^{**})$ 
we have that $x_{1}^{**} \bot y_{1}^{**}$.

Consider now, $\{ \zeta \quad d_{\Gamma}(Y/X)=const \}$ that describes parabolic movements. 
When $d_{\Gamma} \sim \mid \cdot \mid$
we have that $d_{\Gamma}(Y/X)=const$ implies $Y/X \nrightarrow 0$, that is $Y/X =const$ is complementary to 
$ Y/X \rightarrow 0$ and conversely if $Y/X \rightarrow 0$ we have that $d_{\Gamma}(Y/X) \neq const.$
Assume $\{ \gamma_{2} < \lambda \}=M$ with $\gamma_{2}^{\bot}(\gamma_{2})=0$ implies $\gamma_{2} \in M$
where $\gamma_{2} \in (I_{2})$ on a set $V$ and we assume $(I_{2})$ closed with $(I_{2})^{\bot \bot} \simeq (I_{2})$.
If we define $(I_{1}) \subset (I_{2})^{\bot}$ through $d_{1},d_{2}^{\bot}$, we have that $(I_{1})$
can be considered as closed in $(I_{2})^{\bot}$ for instance $d_{2}^{\bot}/d_{1} \rightarrow 0$ with 
$d_{2}^{\bot}(z) \sim d_{2}(\frac{1}{z}) \rightarrow 0$, when $z \rightarrow \infty$. 
Note that if we assume $d_{2}(\gamma) \rightarrow 0$ implies $\gamma \rightarrow \gamma_{0}$ where 
$\gamma$ reduced, we must have $d_{2}(\zeta) \rightarrow 0$ implies $\zeta \rightarrow \zeta_{0}$.

Consider the problem $L(U \gamma,U \gamma)=L(\gamma,\gamma)$ and $L(U \gamma,U^{\bot} \gamma)=0$. An
axes $L(U \gamma,\gamma)=L(\gamma,\gamma)$ corresponds to an axes for $U^{\bot}$ and $L(U^{\bot} \gamma, \gamma)=L(\gamma,\gamma)$.
If the condition $L(U \gamma_{1}, \xi)=0$ implies $\xi=0$, we must have $U=id$, that is $\gamma_{1}$ on the reflection axes.
We assume $\xi=0$ implies $J \xi=0$ since $J$ is assumed to preserve orthogonality, that is 
$L=0 \rightarrow_{J} E=0$. This condition is also assumed for the inverse $J^{-1}$. Thus, if 
$E(V \gamma_{2} , J \xi)=0$ implies $J \xi =0$, we must have $V=id$, that is a reflection axes in $L$,
has a corresponding axes for invariant points in euclidean metrics. In particular 
$H_{1}(\gamma_{2})=< \gamma_{2}, J^{-1} \gamma_{2}>$
where $H_{1}(V \gamma_{2})=< V \gamma_{2}, \gamma_{1}>=< J U \gamma_{1}, \gamma_{1}>=H(U \gamma_{1})$.
$H_{1}(J \gamma_{1})=0$ implies $H(\gamma_{1})=0$ and we can use that $H_{1}(J \gamma_{1})=H(e^{v} \gamma_{1})$
that is $H : (\gamma_{1},L) \rightarrow (\gamma_{2},E)$ (positive) linear functionals.

Given a reflection axes $R$ to $U$, there is a reflection axes $R^{\bot}$ to 
$U^{\bot}$, such that $R \cap R^{\bot} \neq \emptyset$ (note that $R^{\bot}$ is not unique). 
Particularly, if $\gamma(\zeta) \in \Sigma$, where $\Sigma=\{ F(\gamma)(\zeta)=const \}$,
we assume $R^{\bot} \bot R$ with points in common. We have that $R^{\bot}=R^{\bot}(R)$, $R^{\bot} \cap \Sigma = \{ (x_{0},y_{0}) \}$, 
for some $R^{\bot}$. When $F$ is constant on $R$, we have $F \neq const$ on $R^{\bot}$. 
Consider $V,V^{\bot} \in H'$ and use the topological isomorphism $H'(E) \rightarrow Exp(E^{*})$,
where $Exp$ are entire functions or regular with slow growth (\cite{Martineau}). Assume compatibility conditions,
such that $V^{\bot} \gamma_{2}$ regular for $F$ or $F(V^{\bot} \gamma_{2}) = const.$. Then since
$\widehat{V} \cdot 1 = 1 \cdot \widehat{V}$, we have $\big[ V^{\bot},I \big] = \big[ I , V^{\bot} \big]$,
over regular approximations. Further, we assume $JI=IJ$ on the set corresponding to $\Delta$. If 
$R^{\bot}$ is an axis for invariance for $V^{\bot}$, we have an axis $\widehat{R}$ to $\widehat{V}$
and $R^{\bot} \simeq \widehat{R}$, further a line $\tilde{R^{\bot}}$ corresponding to regular approximations.

Assume for parabolic rotation $\tau_{p}$ that $F(\tau_{p} \gamma)=const$, then over $S=\{ F(\gamma)=const \}$, 
$\frac{d}{d \gamma} {}^{t} \tau_{p} \equiv 0$. Note that the condition 
$\frac{\phi_{2}'}{\phi_{1}'}= c \frac{\phi_{2}}{\phi_{1}}$
corresponding to elliptic rotation, that $\frac{\phi_{1}^{2}}{a^{1}} + \frac{\phi_{2}^{2}}{b^{2}}=1$, 
for constants $a,b$. Define $F(x,y)=F_{1}(\frac{y}{x})$ and $\tau_{e}$ elliptic rotation, 
we then have ${}^{t} \tau_{e} F(x,y)={}^{t} \tau_{p} F_{1}(\frac{y}{x})$
and $F_{1}(\frac{y}{x})=F(x (1,\frac{y}{x}))$ or ${}^{t} \tau_{e} F = {}^{t} \tau_{p} F_{1}$.
In particular we consider the domain $D(F)$, that is we can define $D(F)=xD(F_{1})$
$\frac{1}{z} F(\frac{1}{z})=f(z)$ for $0 < \mid z \mid < R$.

Consider $G(z,z^{\bot})=(G(z),G^{\bot}(z))$, where $G(z),G^{\bot}(z)$ are real. Using Tarski-Seidenberg's theorem, if 
$\{ (z,z^{\bot}) \quad G < \lambda \}$ is semi algebraic, we have the same for 
$\{ z \quad G < \lambda \}$ and in the same manner for $z^{\bot}$.
Assume $G G^{\bot}=G^{\bot}G$, then $G^{2} \sim G^{2} - G^{\bot 2} + 2 i G G^{\bot}$. For the real part, 
assume $\{ \frac{G}{G^{\bot}} < \frac{G^{\bot}}{G} \}$ bounded, this corresponds to $G^{\bot 2} \prec \prec G^{2}$
in $\infty$. When the set is unbounded, we have that $G^{\bot} < \lambda$ implies $G < \lambda$ in 
$\infty$. If $\lim_{z \rightarrow z_{0}} G=\lim_{z^{\bot} \rightarrow z_{0}} G$, 
we have that $z_{0}$ is an isolated point. Note that it remains, given $U \gamma$ a movement, 
to define $U \rightarrow \tilde{U}$ and $\gamma(\tilde{U} \zeta)=U \gamma(\zeta)$. 
Consider $G(z) \rightarrow G(z^{\bot})=G^{\bot}(z)$ as a projection operator and assume
$\bot$, dependent of orientation, that is one-sided, then we have that $\{ G < \lambda \}$ semi algebraic 
$\Rightarrow \{ G^{\bot} < \lambda \}$ semi algebraic, can be seen as a transmission property.

Assume $V J = J U$ and $J : p \rightarrow q$. Consider $U \rightarrow V \in (I_{2})'$, 
that is $V$ is a functional over $(I_{2})$. In particular when $( V - I) \gamma=0$ 
and $U \gamma_{1}=\gamma_{1}$ implies $V \gamma_{2}=\gamma_{2}$. Given $(I_{2})$ is closed 
we have that $(I_{2})=({}^{o} I_{2})^{o}$. When $(I_{2}) \subset (I_{He})^{\bot}$ with Runge's property,
we can regard $(I_{2})$ as closed in $(I_{He})^{\bot}$ (with respect to $\dot{B}$).
The moment problem, for $(I_{2})=E_{0}$ provides a continuation to $C$. 
For instance, $\gamma^{\bot}=e^{v} \gamma$ with $e^{v} \rightarrow 0$ (or $\infty$) and 
$F(\gamma,\gamma^{\bot}) \rightarrow 0$
where $\gamma^{\bot}/\gamma \rightarrow 0$ in $\dot{B}$, thus we have existence $F_{1}$, 
such that $F(\gamma,\gamma^{\bot})=F_{1}(e^{v})$.

\section{Intermediate ideals}
Define $V_{1}=\{ x \in \Omega \quad F(\phi x)=F(x) \}$, for a fixed $F$, that is $\phi x - x \in \mbox{ ker }F$, 
an ideal of holomorphy. We assume here $F$ linear in $x$, but it is not assumed linear in $\Omega$. 
Assume $V_{1}$ is given by $\{ f-c_{1} \}$, where $f$ is analytic and $c_{1}$ scalars and in 
the same manner for $V_{2}$. Using the theorem on intermediate values, we can assume existence 
of $V$ between $V_{1}$ and $V_{2}$.
Assume $(J)$ an ideal such that $(J)=\mbox{ ker } \phi_{1}$, for instance $\phi_{1}=\phi-1$ and $(J) \oplus (J)^{\bot}=(E)$, with $\{ 0 \} = N(J) \cap N(J^{\bot})$, 
that is $\phi_{1} + \phi_{1}^{\bot}=id$ (local identity). 
Given that $S_{1} \rightarrow S_{2}$ algebraic, we have type (A) first surfaces, that is we have a 
lifting principle over an algebraic polyhedra (\cite{Nishino75}). 
Concerning existence of $(I_{0})$ such that $(I_{1}) \subset (I_{0}) \subset (I_{2})^{\bot}$. 
Given that the inclusion is continuous (closed) and injective, the existence can be derived 
using the theorem of intermediate values.
A sufficient condition (and necessary) for inclusion between weighted space in $\mathcal{D}_{L^{1}}'$
is that quotient of the corresponding weights goes to $0$ in the $\infty$. If in this setting 
$(I_{1})$ has the weight $\rho_{1}$,
$(I_{0})$ has the weight $\rho_{0}$ and $(I_{2})^{\bot}$ has the weight $\rho_{2}$. 
Then the condition for the inclusion we are looking for is 
$\rho_{2}/\rho_{0} \rightarrow 0$,$\rho_{0}/\rho_{1} \rightarrow 0$ and $\rho_{2}/\rho_{1} \rightarrow 0$.

For the hyperboloid we have given a desingularization $\mu$, that given $F$ meromorphic
then $F \circ \mu$ is holomorphic, we have a lifting principle over a desingularization. 
Consider $(x,y) \rightarrow \frac{y}{x} = \eta(x)$. We denote the diagonal $V=\{ (x,y) \quad y(x) = x \}$.
Given $y=y(x)$ analytic $\in (I)$ and $\eta \in (J)$, we have where $x \neq 0$, that $(J) \subset (I)$ 
and where $x \neq const.$, $V \subset \{ y = const \}^{c}$. $(J)$ is algebraic, if $(I)$ is algebraic.
When $y$ is linear, $V$ is linear.

Assume $x,\eta(x)$ does not have presence of essential singularities in the $\infty$ and consider 
$x(1,\eta) \rightarrow (x,y)$, that is given $\eta$ analytic we have existence $y=y(x)$, 
such that $\eta(x)=y(x)/x$. Given $x$ reduced, we can represent 
a global base.

\newtheorem*{res9}{Intermediate ideals}
\begin{res9}
 Assume $d v$ a measure on the boundary $\Gamma$ that is joint for $(I_{1})$ and $(I_{2})$.
 Define $d \mu_{1}=\rho_{1} d v$ and a corresponding functional $B_{1}(f)=\int f d \mu_{1}$
 and in the same manner for $ d \mu_{2}$ and $B_{2}$, with the condition that $\rho_{2}/\rho_{1} \rightarrow 0$
 when we approach $\Gamma$. If we can find $\rho_{0}$, such that $\rho_{2}/\rho_{0} \rightarrow 0$
 and $\rho_{0}/\rho_{1} \rightarrow 0$ as we approach $\Gamma$, we have existence of a corresponding
 functional $B_{0}$, such that $B_{2} \leq B_{0} \leq B_{1}$ and an ideal between $(I_{1})$ and $(I_{2})$.
\end{res9}

Note that $(I_{Hyp}) \subset (I_{Phe})^{\bot} \subset (I_{He})^{\bot}$. Further if we have existence of
$\alpha$ of bounded variation such that $\int_{\Omega} g d \alpha=0$ and $d \alpha$ reduced,
we have that $\Omega$ is a strict carrier to $\alpha$ (\cite{Martineau}), that is $\alpha$ 
can be represented with compact support.
Assume $d \beta=e^{\phi} d \alpha$, where $\phi \in L^{1}$, here we associate $\phi$ to  a hyperbolic 
movement.
When $\phi$ is linear, we can define the continuation $g \rightarrow \tilde{g}$ using the Fourier-Borel duality.
Note that if $e^{\phi} \rightarrow 0$ on a radius $L$, we have that $e^{\phi} \rightarrow 0$ 
on a disc. In this case $\beta$ can be chosen as summable.
Given existence of $(I_{0})$ such that $(I_{Hyp}) \subset (I_{0}) \subset (I_{He})^{\bot}$ we can chose  $(I_{0})=(I_{Phe})^{\bot}$.

Given a point and a normal in the point, if we consider a neighbourhood of the point with regular boundary, 
then the boundary can be oriented. This corresponds to the concept of a pseudo vector (\cite{Janc89}). 
The regularity conditions for dynamical systems and the corresponding conditions for first surfaces,
for instance the condition (N) (\cite{Nishino68}), can be used instead with advantage. Thus we can replace
the concept of pseudo vectors for first surfaces and transversals, given that we assume the first 
surfaces oriented. Note that this condition is necessary for the transversals to be algebraic.

Note the following problem: Determine a dynamical system corresponding to a hyperbolic system 
such that $\gamma$ in this system has normal $n_{1}$ in $p_{0}$ and with $p_{0}$ a zero to the system. 
We assume further $n_{1}$ locally algebraic. In this case the regularity conditions indicate a pseudo normal orientation 
(\cite{Janc89}). In the same manner for $n_{2}$ and a dynamical system corresponding to a partially hypoelliptic system.
Thus $\Delta \rightarrow n_{j}$ approximates singular points that are given by the right hand side.
When $n_{j}$ are locally algebraic, we have that $\gamma$ is defined as locally analytic.
Assume $d_{j}$ the respective distance functions to the joint boundary (first surfaces) with $d_{2}$
reduced, such that $d_{1}/d_{2}$ a distance.
Using that $d_{1}/d_{2} < \epsilon$, we can use the type (A) condition (\cite{Nishino75}), 
that is assume $S_{1},S_{2}$ first surfaces closely situated, 
we then have that $S_{1}$ intersects the normal in the same manner as $S_{2}$.
In the same manner if $\lambda \bot \eta_{2}$ and $\frac{d \lambda}{d t}=(X_{2},Y_{2})$. 
We consider $(X_{1},Y_{1}) \rightarrow_{\mu} (X_{2},Y_{2})$, as a continuation, 
that is $\mu_{T} \rightarrow \frac{d \lambda}{d T} \mid_{T=0}$ 
when $T \rightarrow 0$ and $\mu_{T} \rightarrow (X_{2},Y_{2})$. 

Assume now $\gamma_{1} \in (I_{1})$ satisfies a dynamical system $(X_{1},Y_{1})$ and in the same 
manner that $\gamma_{2} \in (I_{2})$ satisfies $(X_{2},Y_{2})$. We argue that there exists a dynamical system $(X,Y)$ 
associated to $(I)$. More precisely $(I_{1}) \subset (I_{j}) \rightarrow (I_{2})$ such that $F_{j}$ are 
Hamiltonians corresponding to the approximating ideals $(I_{j})$,
that is $(-\frac{d F_{j}}{d y},\frac{d F_{j}}{d x})$ give the dynamical system we look for in limes. 
We assume for this reason $\eta \bot \frac{d F_{j}}{d x}$ and $\eta \bot \frac{d F_{j}}{d y}$, 
why assuming the conditions in Weyl's lemma, $\eta \in C^{1}$ and $\mbox{ supp } \eta \downarrow \{ 0 \}$. 
In the same manner we can consider $(I) \subset (I_{j}) \rightarrow (I_{2})$ where again $\mbox{ supp }\eta \downarrow \{ 0 \}$. 
In a discussion of desingularization, we note that if $\Omega_{1}$ is associated to $(I_{1})$ 
and $\Omega_{2}$ to $(I_{2})$ we have that we only assume $\gamma_{1} \rightarrow \gamma_{2}$ continuous, 
why we can not establish $\Omega_{1} \rightarrow \Omega_{2}$ as a continuous mapping. Note that a desingularization
is not established for non-hyperbolic operators.

Concerning $MV(x,y)=\frac{x-y}{\log x - \log y}$ (\cite{Vuorinen07}), if $r_{T}'$ is locally algebraic and 
$r_{T}' y(x)=y(r_{T}' x)$, we have $M(r_{T}' x,r_{T}'y)=r_{T}' M(x,y)$. Let 
$MV(x,y)=\frac{x (1 - \frac{d y}{d x})}{-x \frac{d \eta}{d x} \log \eta(x)}$, where 
$\eta(x)=\frac{y(x)}{x}$ and $=\frac{-y (1 - \frac{d y}{d x})}{x \frac{d \eta}{d x}}$, 
without degenerate points in the denominator. Given that $g=\frac{d}{d x} e^{\phi}$, we then have
$\frac{G(g)}{g} \leq \frac{M(g)}{g} \sim \frac{d x}{d \phi}$, that is we have a lower bound 
for the inverse mapping to the phase. If we assume $f_{j}=e^{\phi_{j}}$ with $\phi_{j} \in L^{1}$ for $j=1,2$ and 
$\frac{\phi_{1}}{\phi_{2}} \rightarrow 0$, we have $\phi=\phi_{1} + \phi_{2}$ and $\phi_{1} \bot \phi_{2}$, 
which implies $e^{\phi}=f_{1}f_{2}$. 
In (\cite{Riesz15}) we note Theorem 17, that is if a Dirichlet series 
is summable through arithmetic mean, it is summable through logarithmic mean.
The proof is based on a discussion on the function $(\frac{\log w - \log t}{w-t})^{k-1}$ 
$=(\frac{1}{L(w,t)})^{k-1}$ that is increasing steadily from $t=1$ to $t=w$ and the limit 
as $t \rightarrow \infty$ is $w^{1-k}$. Riesz gives a proof for $k>0$. For the converse result 
(in general false), we refer to theorems 19,20 (\cite{Riesz15}).

Consider now $\gamma \rightarrow \eta=y/x$, then we have that reflection through $\gamma \rightarrow -i \gamma$
can be written $\eta \rightarrow 1/\eta$ and in the same manner reflection through the real axes, $\eta \rightarrow - \eta$.
The condition on vanishing flux is then $\int_{\beta} d F(\frac{1}{\eta})=0$. Assume $\eta_{2}$ denotes 
continuation of $\eta_{1}$ to $(I_{1})^{\bot}$, such that $\eta_{1}^{*} \sim 1/\eta_{2}$. Assume 
$F(\eta_{1}^{*})=F^{*}(\eta_{1})$, we then have if $F(\eta_{1})-F^{*}(\eta_{1}) \sim F(\eta_{1}) - F(1/\eta_{2}) \sim 0$, 
that $F$ is symmetric over the continuation and we have a transmission property in this case.

Assume now $G(f)=e^{M(\phi)}$ for $f=e^{\phi}$ and $g=\frac{d}{d x} f=e^{\psi} =\frac{d \phi}{d x} M(g)$
that is $\frac{M(g)}{g} \sim \frac{d x}{d \phi}$. Assume $\phi,\psi$ sub-harmonic. Given a lower bound for 
$\frac{d x}{d \phi}$, we can conclude that if $\phi \rightarrow \phi_{0}$, we see that $x \rightarrow x_{0}$. 
Given $1 \leq \frac{d x}{d \phi}$, we have that $\phi \leq \phi \frac{d x}{d \phi}=z(\frac{1}{\phi})$, 
given $x$ is analytic considered as a function of $\phi$ in $0 < \mid \phi \mid < \infty$. Given 
$z (\frac{1}{\phi}) \rightarrow 0$, we have that $\phi \rightarrow 0$, why we have existence of $w$, 
such that $\frac{1}{w}(\frac{1}{\phi}) \rightarrow 0$ implies $\phi \rightarrow 0$, that is $\frac{1}{w}(\phi) \rightarrow 0$ 
implies $\phi \rightarrow \infty$, why $w$ is downward bounded. Thus we have where $x$ is analytic, 
a mapping $x \rightarrow w$, such that $w$ is locally 1-1, that is if $w \rightarrow w_{0}$ implies 
$\phi \rightarrow \phi_{0}$, it is sufficient to consider the phase. Note that 
if $\frac{d \phi}{d x} \rightarrow 0$ implies $x \rightarrow 0$, that is
$f'/f \rightarrow 0$ implies $x \rightarrow 0$, then we have that $f \in (I_{He})^{\bot}$ for 
all $x \neq 0$.

\section{The transmission property}

\subsection{The transmission property for summable distributions}
Immediately, in the case when
$F$ can be represented by a linear functional, we can define $\bot$ using for instance annihilators. 
If $I(f)=\delta_{0}*f$ we have that $I$ reflects the support. If $f$ has symmetric support we have that 
$I(f) \sim f$. If further $f$ is proper, then $I(f)$ has the same property.
Assume in particular $U_{j}$ closed in some set $Z$, assuming Schwartz type topology, 
we have that $(x,\varphi) \in V \times U_{j}$ is closed iff $x \rightarrow \varphi(x)$ continuous. 
If singularities are approximated through transversals $L$, there are points in $Z \cap L$ 
that do not belong to $U_{j}$. More precisely, if $u_{0} \in L \cap U_{j}$ and 
$u_{j} \rightarrow u_{0}$, we have that $u_{j} \in Z \backslash U_{j}$.
Define $F(I x)=F(-x)$ and ${}^{t} I F(x)=-F(x)$, this can be seen as a weak definition of odd 
operators. Thus if $F I = I F$ we have that $F$ is odd. If $F I = F$ we have that $F$ is even.
Note that $\mbox{ supp } \big[ I,E \big]=- \mbox{ supp }E$. If $E,F$ are both one-sided, 
we may have that $\big[ I,E \big] \big[ I,F \big]$ is one-sided
but $\big[ I,E \big] F \equiv 0$. 

In the terminology of oscillating integrals, $F(x)=\int a(x,\theta)e^{i s(x,\theta)} d \theta$ 
where $a(x,\theta)$ are poly homogeneous, for $x \in X$ and $\theta \in \mathbf{R}^{N}$. 
The phase function $s(x,\theta)$ is assumed real, 1-homogeneous in $\theta$ with 
$s_{m}(x,\theta) \neq 0$, as $s \neq 0$. Assume $S=\{ s_{\theta}(x,\theta)=0 \}$ closed and conical 
and $WF(u)=\{ (x,\theta) \quad \theta=s_{x} \quad s_{\theta}=0 \}$. Using a  bijection 
$\varphi$, $WF(u (\varphi(x)))=\{ (\varphi(x),{}^{t} \varphi'(x) \theta) \in WF(u) \}$. Assume 
$x^{*}=\xi$, $\eta^{*} \sim \varphi(x)$, why $\frac{d}{d x} s(\varphi(x),\xi)=\frac{d \varphi}{d x} s_{x}=\varphi'(x) \xi$, 
then we have $\eta = \frac{d \varphi}{d x} \xi$ and ${}^{t} \big[ \frac{d \varphi}{d x} \big]^{-1} \eta=\xi$. 
In particular $\int a(x,\xi) e^{i < x,\xi> - H(\xi)} d \xi \sim \widehat{(A/H)}$ (\cite{Garding87}).

Note for the dual $a'$ to a polynomial, if $a'-a \equiv 0$ on the radius in a disc, 
we have that $a'-a \equiv 0$ on the disc. 
For a pseudo differential symbol $a$ we must except the imaginary axes. 
Assume that $a$ has a symmetric kernel with respect to the imaginary axes, for instance
$a'-a = \Sigma_{j} \delta^{(j)}$, where $\delta$ is assumed to have support on
the imaginary axes. 
Given a very regular boundary and $\eta_{2}$ reduced
we can consider $1/\eta_{2} \in \dot{B}$, why $(I_{1})$ can be seen in $\mathcal{D}_{L^{1}}^{\bot}$ 
that is $P(D) \tilde{f}_{0}$
is associated to $1/\eta_{2}$, which can be seen as a topological transmissions property.

\subsection{The transmission property for two bases}
The transmission property, when it is derived from involutive reflections $U$, implies 
presence of a normal model. In the last section, we see that a global model does not imply a normal
model. If we can map bijectively a broken ray
onto a transversal, we have one-sidedness for one segment with respect to the first ($x-$) axes and one-sidedness 
for another segment with respect to the second ($y-$) axes.
Now consider one-sidedness with respect to two planes, assume $E(x,y)$ the kernel corresponding to 
the first segment and $F(y,z)$ the kernel to the last segment. Assume that we have existence of 
a segment $ y \rightarrow y'$
(transversal) according $\frac{d y'}{d y}=1$. We assume again the segment normalised using distances.
We write $G(y,y')$ for the kernel corresponding to the middle segment. 
In order to compare with one-sidedness, we assume $G(y,y') \rightarrow 1$
as $y \rightarrow y'$. Note that $F$ is assumed symmetric after reverting the orientation 
for the last segment.
If the composition is algebraic, we do not have that the factors are algebraic.

Assume we use distances to represent the convergence, let $d_{1}$ be the distance 
to $\xi_{2}=0$ and $d_{2}$ the distance to $\xi_{1}=0$.  When the distances are reduced we can assume the inverse 
a distance, otherwise a pseudo distance. Assume $A \rightarrow A' \rightarrow B' \rightarrow B$, where $A,B$ on $\mid \xi \mid=1$ and
$A',B'$ on axes. Given symmetry for the kernel, we can assume $d_{2} \rightarrow 0$ has a 
reflection in $d_{1}^{-1} \rightarrow 0$.
If $d_{1}(\xi) \rightarrow 0$ as $\xi \rightarrow 0$, then $\frac{1}{d_{1}}(\frac{1}{\xi}) \rightarrow 0$, 
when $\xi \rightarrow 0$
that is $\overline{\xi}=1/\xi$. 

Let $i(d_{j})=1/d_{j}$ that maps distances onto pseudo distances. 
If $A \sim B$
(congruence) we can assume $d_{1}d_{2} \sim 1$. In particular $d_{1}=0$ for $z=A$,
$\frac{1}{d_{1}}=d_{0}=0$ for $z=A'$, $\frac{1}{d_{0}}=d_{2}=0$ for $z=B'$ and $\frac{1}{d_{2}}=0$, 
for $z=B$.
In this case we have that $d_{1}d_{2}=0$ implies $z =A$ or $z=B'$.
As $\frac{d_{1}}{d_{2}}=0$ we have $z=A$ or $z=B$. Under the symmetry condition for $f_{\lambda}=f-\lambda$ we have
$\{ d_{1} \mid f_{\lambda} \mid < c \}$ $\sim \{ \frac{1}{d_{2}} \mid f_{\lambda} \mid < c \}$, that is we can assume 
$f_{\lambda}$ has algebraic growth at the boundary.
If we instead assume $A$ singular in the sense that $ d_{1} \frac{1}{\mid f_{\lambda} \mid} < c$ at the boundary, 
we have that $f$ is downward bounded at the boundary and $\{ \mid f_{\lambda} \mid < c \}$ can be seen 
as bounded. Simultaneously if $B$ is regular, we have that the corresponding set is unbounded at 
the boundary.

Consider $f$ with algebraic continuation in $\{ d_{1} \frac{1}{d_{2}} \mid f \mid < c \}$. We are assuming
$\tilde{f}(z,w)$, where $z$ is the reflection of $w$. 
Assume now $z,w$ a continuation of path $A \rightarrow A' \rightarrow B' \rightarrow B$. 
Define instead $d_{1} \rightarrow 0$ as $z \rightarrow A$ and $d_{1} \rightarrow 1$ as $z \rightarrow A'$
and $\frac{1}{d_{2}} \rightarrow 0$ as $z \rightarrow B$ and $d_{2} \rightarrow 1$ as $z \rightarrow B'$.
We assume $d_{0}(A')=d_{0}(B')=1$, we have that $d_{0} \frac{d_{1}}{d_{2}} \mid f \mid(z,w) \rightarrow $
$0$ as $z \rightarrow A$, $1$ as $z \rightarrow A'$, $1$ as $w \rightarrow B'$ and $0$ as $w \rightarrow B$.

Starting from a point at the boundary $A$, assume that we have existence of a path $\tau$ between 
$A$ and $B$, that does not contain any other singularities than $A$, in the sense that $A$ can 
be identified with a regular point $B$, we could represent $A$ and $B$ in the same leaf. 
The points can be seen as isolated on the path. If not all points on the boundary
$\Gamma$ are singular, we can always find a path to $B$, if all singularities are situated on the 
boundary, we have that the path is regular.

\subsection{General remarks}
The condition $\log f \in L^{1}$ can be seen as $\log \mid f \mid \in L^{1}$, that is 
$\mid f \mid = \mid f^{*} \mid$ is a symmetric condition. For selfadjoint operators, it is sufficient
to conclude hypoellipticity, to consider translation invariant sets.

Polynomial operators have solutions with one-sided support, that is we can assume $F(\gamma)(\zeta)=I(\zeta)$, 
where $F=\widehat{E}$ and $E$ has one-sided support, when
$\gamma$ is polynomial outside the kernel to $E$. Thus we have for the continuation 
$\gamma_{1} \rightarrow \gamma_{\delta}$ algebraic
that the corresponding $F_{\delta}$ can be selected with one-sided support. 
Starting with $F(e^{-v} \gamma_{1})=F(\gamma_{2})$, with $v \in L^{1}$, such that $v \equiv 0$ over the
set $\Delta$ where $\gamma_{1}=\gamma_{2}$. Assume $V \gamma_{2} \sim \eta$ and $\eta \rightarrow \eta^{*}$
reflection through the axes $\eta=\eta^{*}$, $\eta = y/x$. Assume $v(\eta,\eta^{*}) \in L^{1}$ and note
that $v=0$ implies $\mid \eta \mid \leq 1$ ($\geq 1$). Let $H(\alpha)=\int H(\eta,\eta^{*}) \alpha(\eta) d \eta$,
then for $\alpha \in \dot{B}$, we can assume $H \in \mathcal{D}_{L^{1}}'$ a Schwartz kernel. The condition
$v \in L^{1}$ means isolated singularities and we can assume $v \sim v_{1}$, where $v_{1}$ is analytic.
When $H \sim_{m} H_{1}$ locally algebraic, we do not have a normal model for $H_{1}$, even when $\eta \rightarrow \eta^{*}$
is involutive. However we have a global model according to the moment problem.

\newtheorem*{res12}{A global representation}
\begin{res12}
 The representation $F(\gamma_{1} \rightarrow \gamma_{2})(\zeta) = f(\zeta)$ and $F(\gamma_{1})=F(e^{v} \gamma_{2})$,
 with the condition that $v \in L^{1}(\eta,\eta^{*})$, gives a global model. It is not necessarily a 
 normal model.
 \end{res12}

Note that $v$ is defined by the movement.  By considering $w=v + v^{\bot}$, where $v^{\bot}$ is defined
by the orthogonal movement and where both $v,v^{\bot}$ are assumed continuous, if we assume the 
compatibility condition $w=0$ on invariant points and on the 
intersection of the supports to $v,v^{\bot}$, we can obviously give a global representation for $w$.
The conclusion is that if the singularities are reached using $w$, we have a global base for the  model.
 
Assume $A=A^{*}$ a hypoelliptic ps.d.o, then $A$ has very regular parametrices, with trivial kernel.
Thus the parametrices to $A$ have a transmission property modulo $C^{\infty}$, in the sense that the regularity
is symmetric with respect to diagonal. For a partially 
hypoelliptic operator, we do not have a trivial kernel, but the Schwartz kernel has hypoelliptic 
representation over kernel (zero space) and outside
the zero space (\cite{Dahn13}). Obviously we do not have that the sum of kernel has the symmetric regularity property 

Note that if we define $\bot$ using $g(\frac{1}{x})=0$ on $V \ni 0$ (a bounded set) 
we must separate between the case where $g \sim \frac{\mbox{ Im} f}{\mbox{ Re }f}$ analytic 
and the case where $g \equiv 0$ in a disc-neighbourhood of $0$. When $\bot$ is defined by 
a rational function $g$ we have that $\bot$ is of order 0, Further, when $g$ is analytic we have that 
$\bot$ is Cartier.


\subsection{Final remarks}

Starting with the observation, that if $\varphi$ is a closed form in the plane and $\varphi^{*}$ (harmonic conjugate)
with $\varphi^{*}=-i \varphi$, we have that $\varphi$ analytic. A mapping that maps 
$(x,y) \rightarrow (y,-x)$ is pure, that is preserves analyticity. 
Assume $L$ the line $x=y$ (the light cone), then reflection through $L$ can be seen as parabolic 
movements. Assume $L_{0}$ the positive
imaginary axes, $L_{1}$ line $x=y$, $x>0$, $L_{2}$ the positive real axes, 
$L_{3}$ the line $x=-y$ and so on.
Assume $c_{0}(x,y)=(-x,y)$ that is reflection through $L_{0}$ and so on. We then have that 
$c_{j+1}c_{j}(x,y)=(y,-x)$ that is pure mappings.
Further, $c_{j}c_{j+1}(x,y)=-(y,-x)$ and so on.

Consider $\frac{f_{2}}{f_{1}}=\frac{1}{\varphi}$. Note that given $W$ planar in $\mathcal{O}_{AD}$
(\cite{AhlforsSario60}), that is we have existence of $u$ analytic with finite Dirichlet integral,
such that $u$ is constant on $W$.
We then have for every single valued function $u$, that $u$ is linear. We can in this context consider
$\frac{1}{\varphi}$ as a function of $f$. For instance $\frac{1}{\varphi}(f+g) \rightarrow 0$ 
in $\infty$, when $f_{1} \bot f_{2}$ and $g_{1} \bot g_{2}$. 

When a domain is defined by a distance $d$, the domain is said to be $d-$ pseudo convex if 
$- \log d$ is plurisubharmonic.
Note that $d_{j}(x) - d_{j}(y)=0$ iff $e^{d_{j}(x)}-e^{d_{j}(y)}=0$, $\mid d(x) - d(y) \mid \leq d(x-y)$.
In particular, on the set where $d^{*}=d$ and $d^{2} \sim 1$ we have if $d (f)=0$ that 
$\mid f \mid =1$ and $\frac{d}{d x} d(f)=0$ implies $f=1$. Assume $d>0$ a distance on the phase space 
to symbols, then we have
$d(\phi)=d(\psi)$ implies $\tilde{d}(e^{\phi})=\tilde{d}(e^{\psi})$, where $\tilde{d}(e^{\phi})=e^{d(\phi)}$.
When we have equality in the tangent space, we have $e^{\phi}=C e^{\psi}$. Further if $\tilde{d}$ 
is locally 1-1,
we have that $\tilde{d}>0$ where $d>0$. When $d<0$ we can consider 
$\tilde{d}(f) \tilde{d}(\frac{1}{f})=e^{d(f) - d(f)}=1$, that is $\tilde{d}(\frac{1}{f}) \sim 1/\tilde{d}(f)$ 
and $\tilde{d}$ can be seen as ''algebraic``. Note that this means that $d$ corresponds to $(f_{1},f_{2})$, $\overline{d}$ to $(f_{1},e_{2})$, 
$-d$ to $(e_{1},e_{2})$,
$-\overline{d}$ to $(e_{1},f_{2})$, $i d$ to $(e_{2},f_{1})$, $i \overline{d}$ to 
$(f_{2},f_{1})$, $-i d$ to $(f_{2},e_{1})$ and $- i \overline{d}$ to $(e_{2},e_{1})$. 
We have in this case a two-sided solution $f_{1}e_{1}=e_{1}f_{1}=1$. 

Assume $A \rightarrow B$ through reflection points $s_{1},s_{2}$. Assume 
$\frac{1}{2} d \overline{d}$ $=\frac{1}{2} (d_{1}^{2} + d_{2}^{2})$ defines the distance between 
$A$ and $B$. Given that the domain for $d \overline{d}$
is connected and $d \overline{d}$ monotonous on the path, there is a $\mu$ between $s_{1}$ and 
$s_{2}$, using the theorem of intermediate values. For instance if $d_{1}=1$ in $s_{1}$ and 
$d_{2}=1$ in $s_{2}$, we have that $\frac{1}{2} d \overline{d}=1$ on the  path between $s_{1}$ and 
$s_{2}$. If $d_{1}^{2}$,$d_{2}^{2}$
are polynomial, we can consider a semi algebraic set and the path $\mu$ between the reflection 
points as a part of the boundary.

\section{The boundary}

\subsection{First surfaces}

First surfaces are naturally invariant for all movements and existence of a regular approximation 
outside the first surface implies existence of a $V$ such that $\frac{d V}{d T} \neq 0$ over 
this set.
More precisely, assume the first surfaces are defined by movements with the compatibility condition 
$U^{\bot} \gamma \subset S_{0}$ or $U^{\bot} \gamma \cap S_{0} = \emptyset$, where $S_{0} = S \backslash \{ (x_{0},y_{0}) \}$, and $(x_{0},y_{0})$ is a point
that is invariant for both $U$ and $U^{\bot}$. When $V$ is translation we have that $V^{\bot}$ 
can be seen as a transversal,
when $V$ is an elliptic movement, then $V^{\bot}$ is a line through $0$, when $V$ is scaling 
then $V^{\bot}$ is scaling.
Note that the  invariance principle does not include the movements $\bot$ lorentz movements. 
However we can define $U^{\bot} \rightarrow V^{\bot}$ using Radon-Nikodym's theorem. 
Note that $L(U \gamma, \xi)=0$ iff $L(\gamma,{}^{t} U \xi)=0$. If $\xi=U^{\bot} \gamma$
we see that ${}^{t} U=(U^{\bot})^{-1}$ implies $\gamma \in K$ (the light cone).

Starting with a normal tube, where $L$ is an analytic line ($\Delta$), that is assumed to intersect 
the first surfaces transversally, we apply a movement $U \rightarrow V$, with the property that $\gamma \in \Sigma$ 
implies $V \gamma \subset \Sigma$, where $\Sigma=\{ F(\gamma)=const \}$. We consider a representation 
$F \sim G$, that is regular and with equivalent first surfaces. We have that (\cite{Nishino75}) 
the bases $\gamma$ can be selected as algebraic for $G$, given that $G$ has first surfaces of type (A). 
Thus, we can regard $F(V^{\bot} \gamma)$ as an associated function to $F(\gamma)$ (and $F(U \gamma)$). 
Assume, $V$ an euclidean movement, that preserves first surfaces to $f$ and $e^{v}F$ the corresponding associated function, 
then we have that $e^{v}F$ has type (A), if $f$ has type (A).
More precisely given $L=\{ \gamma_{1}=\gamma_{2} \}$ (formed on $\Delta$) analytic, we can define
the boundary $\Gamma$, as first surfaces that intersect $L$ transversally. We can continue $\Gamma$ 
with first surfaces $\Sigma \ni \gamma \rightarrow V \gamma \in \Sigma$. Given the compatibility 
conditions, we can further continue the boundary by replacing $L$ with $V^{\bot}$ and repeat the 
argument.

A simple very regular boundary is, $F_{T}$ holomorphic in $T$ and $\frac{d}{d T} F_{T}$ holomorphic  in $T$, 
where $T \notin \Sigma$ and $\Sigma$ are isolated points close to the boundary. Regular approximations 
are formed as $V_{1} \cup V_{2}$, where $V_{1}=\{ F_{T} \quad \mbox{not const} \}$ and 
$V_{2}=\{ \frac{d F_{T}}{d T} \quad \mbox{ not const} \}$.

Note that given $f$ holomorphic with first surfaces $\{ S_{j} \}$ and $S_{0}=\{ f - \alpha_{0} \}$
and let $E = \lim_{j} S_{j}$, with $E \cap S_{0} \neq \emptyset$, this implies $S_{0} \subset E$. 
The type for the class of first surfaces depends on the representation of the symbol. For instance 
$F(\gamma)=\int A d \gamma$ is not dependent of $V \gamma$, where $V$ is translation. When $F(\gamma)=G(\frac{y}{x})$, 
where $G$ is assumed regular, the representation is not dependent of $V \gamma$, when $V$ denotes scaling. Finally 
$F(\gamma)=F(\mid \gamma \mid)$ is not dependent of $V \gamma$, where $V$ is rotation and when the 
singularities are given by the condition $\log \mid f \mid \in L^{1}$
they are not dependent of rotation of the symbol.

\newtheorem{res13}{Boundary condition}
\begin{res13}
 We assume all singularities are situated at first surfaces for $f(\zeta)$. We assume $\Delta$ corresponds
 to a set $\Sigma$, where $\gamma_{1}=\gamma_{2}$. The condition $\log \mid f \mid \in L^{1}$ implies algebraic singularities.
 When the boundary is represented in $\mathcal{D}_{L^{1}}'$, through $B(\gamma)=\int \gamma d \mu$, where
 $d \mu$ is assumed very regular, if we let $F(\gamma_{1})=B_{1}(\gamma_{1})$, we assume as involution
 condition $d \mu_{1}/d \mu_{2} = e^{v}$, where $v \in L^{1}$ and $v=0$ over $\Sigma$. Further, we assume
 $d \mu_{1}=0$ iff $d \mu_{2}=0$, which when $d \mu_{j}=X_{j} d x_{j} + Y_{j} d y_{j}$, $j=1,2$, corresponds to
 Lie's point continuation.
\end{res13}

\subsection{Radon-Nikodym's theorem}
Starting with the lineality we can define the boundary as follows. Assume $\Delta^{\bot}$ are 
orthogonals (tangents to $\Gamma$) to $\Delta$, we can describe $\Gamma$ using 
$e^{\varphi} \gamma$ regular relative $\Delta^{\bot}$. 

The Radon-Nikodym's model is, given that $I(\frac{y}{x}) \rightarrow 0$ implies $A_{1}(\frac{y}{x}) \rightarrow 0$, 
we have $A_{1}(\frac{y}{x})=I(\alpha \frac{y}{x})$ for some $\alpha \in L^{1}$.
Note that as $A(0,0)=(0,0)$, we must consider $A_{1}(\frac{0}{0})=\frac{0}{0}$, (compare L'Hopital).
Assume $F_{1}(\eta(x))$ where $\eta(x)=y(x)/x$, we then have that 
if the limits are on the form $F_{0/0}$, we must calculate $F_{1}(\frac{dy}{dx})$, 
that is we assume that we have non degenerate points and $x \frac{d \eta}{d x} \neq 0$.

Consider $(I_{1}) \subset (I_{0}) \subset (I_{2})^{\bot}$ and $\gamma_{1} \rightarrow \widetilde{\gamma_{1}} \rightarrow \gamma_{2}^{\bot}$
that is $< \gamma_{2},\gamma_{2}^{\bot} > \sim 0$ (or $\sim 1$), if $\gamma_{2}$ polynomial, 
we do not necessarily have $\gamma_{2}^{\bot}$ polynomial, but it can always be represented as 
an annihilator (or $(1-\gamma_{2}^{\bot})$). The moment problem means,
given $\int g d \mu=0$ with $g \in (I_{1})$ analytic, that the solution if it exists 
is in $(I_{2})$, for instance $d \mu=P(D) d v$, with $d v$ reduced implies $d \mu \in \mathcal{E}'$. 
If $W$ is open in $V$ we have that $H(V)$ is dense in $H(W)$, then $W$ has Runge's property in 
$V$. 
Given Runge's property we can select
$(I_{0})$ in $C_{0}$, that is $d \mu \in C_{0}'$, a Radon measure. Oka's property, in this context
existence of a global base,
is dependent on the domain. 

Consider the problem of localisation: assume that $V$ a geometric set and $f(\zeta)$ 
a symbol over $\Omega$ and $V$ not a subset of $\Omega$. Further that we have $F(\gamma)(\zeta)=f(\zeta)$ 
over $\Omega$. We can construct $\gamma_{1}$ over $\Omega$,
such that there exist $\gamma \rightarrow \gamma_{1}$ continuous (continuation), that is 
$F(\gamma) \rightarrow \gamma \rightarrow \gamma_{1} \rightarrow V$ continuous. Assume for instance 
that $\gamma_{1}^{\bot} \rightarrow V$ continuous (closed) and locally 1-1,
such that $\gamma^{\bot}=r_{T}' \gamma_{1}^{\bot}$ surjective, existence of $\gamma_{1}^{\bot}$
defines $\gamma_{1}$ (as an annihilator), why $\gamma_{1} \in (I_{He})^{\bot}$ implies 
$\gamma_{1} \in C_{0}$, using the moment problem.

The boundary is given by first surfaces $\{ \gamma \quad F(\gamma)(\zeta)=const \}$, 
dependent on the involution condition and by $\{ \zeta \quad F(\gamma)(\zeta)=const \}$. 
When $\gamma$ is algebraic, we can assume a continuous mapping between these sets. 
Let $\big[ F , E_{\delta} \big](\gamma_{2})=const$, as $\delta \rightarrow 0$
and $\log E_{\delta} \in L^{1}$. Then $F$ corresponds to a solution to a (partially) hypoelliptic 
operator and we see that $\mbox{ ker }E_{\delta}$ is the kernel to $J \gamma_{1}$. 
We thus assume ${}^{t} J_{\delta} F(\gamma_{1})=
F(J_{\delta} \gamma_{1})=\big[ F,E_{\delta} \big](\gamma_{2})$. Assume $F$ such that 
$\mbox{ ker }F=\{ 0 \}$ which implies $\mbox{ ker } F \subset \mbox{ ker } {}^{t} J_{\delta} F$. 
Further, $E_{\delta} \rightarrow 1$, as $\delta \rightarrow 0$.
Assume first surfaces to $F(\gamma_{1})$ 
corresponds to a desingularization, then we have $\tilde{S}_{j}=\{ \zeta \quad \big[ F,E_{\delta} \big](\gamma_{2})=const \}$
assuming ${}^{t} J_{\delta} : const \rightarrow const.$

The condition for the normal $n$, $d t=\sigma d n$ can be compared with Radon-Nikodym's theorem. 
We assume in this case $\sigma \in L^{1}$ over an unbounded set. When $\sigma$ is independent 
of some variable, that is we have a planar obstacle, the condition is not satisfied.
When $\sigma$ have bounded sublevel-surfaces, it can be given by a regular function.

\subsection{Distance functions}
Assume $d_{T}(\gamma)=0$ implies $r_{T}' \gamma \in \Delta$ and that $\Delta'$ defined through the symmetry
condition $F(r_{T}' x,y)=F(x,r_{T}' y)$, that is we assume $\Delta \subset \Delta'$. We define 
$d_{T}'(\gamma)=0$ implies $r_{T}' \gamma \in \Delta'$. Define $\Delta''$ through $F(r_{T}'x,y)=F(-y,r_{T}'x)$, 
that is invariance for harmonic conjugation. We define the corresponding pseudo distance $d''$. 
Finally, we define $d^{0}$ as a pseudo distance to $K=\{ x=y \}$. Thus $(r_{T}' x,y) \rightarrow (x,r_{T}' y)$, 
that is reflection through $K$ can be seen as $d^{0}=const$. Further, $(r_{T}'x,y) \rightarrow (-y,r_{T}' x)$ 
can be seen as $d''=const$, $(-y,r_{T}'x) \rightarrow (-r_{T}'x,y)$ through $d^{0}=const$. The composition, 
that is $d'' d^{0}=const$ contains $(r_{T}' x,y) \rightarrow (-r_{T}' x,y)$. 

Note if $d_{1}d_{2}=const$ and $d_{1},d_{2}$ distances, we have that $d_{1},d_{2}>0$ and 
$d_{1}d_{2}=0$ implies that $d_{1}=0$ or $d_{2}=0$, thus $A \sim B$ (congruence), that is 
$d_{1}$ denotes the distance to $A$ and $d_{2}$ denotes the distance to $B$ we have that 
$z=A$ or $z=B$. In particular if $A=0$ and $B=\infty$, we have that 
$\frac{d_{1}}{d_{2}}=0$, for $\frac{1}{d_{2}}(z)=d_{2}(\frac{1}{z})$.
The example $d_{1}d_{2}=const$ with $d_{1} \rightarrow 0$ as $z \rightarrow \Gamma$ under 
the condition $d_{2}(\frac{1}{z})=1/d_{2}(z)$ implies $1/d_{2}(z) \rightarrow 0$ and if 
$d_{2}$ reduced, we have $\frac{1}{z} \rightarrow \infty$, that is $z \rightarrow 0$. 

We are assuming no essential singularities in the $\infty$, why if $1/d_{\Gamma} \rightarrow 0$, 
as $z \rightarrow \infty$, it is a boundary point but not a singularity. 
Note that the condition $\log f \in L^{1}$, means that for the ideals defined by $d_{j}^{2}$
polynomials, such that $d_{j} \rightarrow d_{j+1}$ compact, the singularities in the zero space
to $\lim d_{j}^{2}$, can be algebraic.

Assume $V=\{ \frac{1}{v}=const \}$, we then have that $V$ is bounded if 
$\mid \xi \mid^{\sigma} \leq C \frac{1}{v}$ when
$\mid \xi \mid \rightarrow \infty$, that is if $v$ algebraic in $\infty$ we have that 
$V$ is regular and bounded in $\infty$. If $V' = \{ v=const \}$ and $v$ 
downward bounded in $\infty$ and in the same way if $\frac{1}{v}$ is algebraic in 
$\infty$ we have that $V'$ is regular and bounded.

Assume $F(\gamma) \rightarrow \gamma_{1} \rightarrow \gamma_{2} \rightarrow \zeta$, 
where $\gamma_{1}$ is hyperbolic and $\gamma_{2}$ is partially hypoelliptic. The problem is to 
determine if $\gamma_{1}$ and $\gamma_{2}$ can be constructed starting from the same lineality $\Delta$. 
An important difference between hyperbolic and partially hypoelliptic symbols, is that $\gamma^{N}_{1} \in (I_{1})$ 
iff $\gamma_{1} \in (I_{1})(\Delta)$ and $\gamma^{N}_{2} \in (I_{2})$ does not imply that $\gamma_{2} \in (I_{2})$. 
Analogously, we have that $d \gamma_{1} \in (I)$ implies $\gamma_{1} \in (I)(\Delta)$ and $d \gamma_{2} \in (I)$ 
is implied by $\gamma_{2} \in (I)$, but not conversely, that is for partially hypoelliptic operators, 
we have $\Delta_{(j)} \downarrow \{ 0 \}$ and for hypoelliptic operators, we have that $\Delta_{(j)}=\Delta$. 
The proof of this is based on the condition that $\Delta \subset$, a domain of holomophy.

Assume $(I_{1}) \subset (I_{He})^{\bot}$ is defined by distances, $d_{2}^{\bot}/d_{1} \rightarrow 0$, 
for $\zeta \in \Sigma$ and $d_{1}=\inf_{\zeta} \mid f-c \mid$, that is distance to the first surface 
$f(\zeta)=F(\gamma)(\zeta)$. Then if $f$ is an entire function, we have that $f \rightarrow c$, as $d_{1} \rightarrow 0$. 
Note that as $\gamma$ is a pseudo differential operator, we can have that $\{ \gamma < \lambda \}$ 
is not relative compact. Further, if $d_{2} /d_{2}^{\bot} \rightarrow 0$, then $d_{2} / d_{1} \rightarrow 0$ 
(one-sided orthogonality). 

Note that $H(U \gamma) \equiv 0$ defines a regular domain, if $H$ has an analytic kernel, that is
$H (U \gamma) = \int U \gamma d \mu =0$, for some $d \mu$ defined by the boundary. Isolated singularities
implies that $U$ is monotropical with an analytic $H$. 
Assume that all singularities can be approximated by
$U^{\bot}$, we then have that given a fixed $U$, we can assume $U \gamma_{1},V \gamma_{2}$ algebraic.
If for $F \in \mathcal{D}_{L^{1}}'$, $F(\gamma_{1} \rightarrow \gamma_{2})$ is analytic over $(I_{1}) \rightarrow (I_{2})$, we assume 
merely continuity for $F$ over $U^{\bot} \gamma_{1}$. Given $\gamma_{1},\gamma_{2}$ a polynomial, 
we can chose $F$ analytic over $\gamma_{1} \rightarrow \gamma_{2}$.

\subsection{Localization}
If the boundary is given locally by one single function, we can represent this function as the 
solution to a differential operator (boundary operator). 

Assume $f$ analytic on a domain $\Omega$, assume $\gamma_{1} \in (I_{Hyp})(\Omega)$ and 
that we have existence of $F_{1} \in H'(\Omega)$ such that $F_{1}(\gamma_{1})=\gamma_{2} \in (I_{Phe})$. 
Given that we have existence of $F_{2}$ analytic over $\gamma_{2}$,
such that $\big[ F_{2},F_{1} \big]$ analytic over $\gamma_{1}$ we have a lifting principle 
over $\gamma_{1} \rightarrow \gamma_{2}$. The base is reversible, if $\big[ F_{1},F_{2} \big]$ 
analytic over $\gamma_{2}$. 

For hyperbolic spaces we have that every pair of points can be linked through a chain of 
analytic discs (the transversal analytic). Assume $S_{j}$ a first surface to a hyperbolic 
base and $\tilde{S}_{j}$ the continuation to the partially hypoelliptic base
as simply connected domains. If we assume the continuation algebraic, with 
$\log \tilde{f} \in L^{1}$, we still have singularities of finite order, 
but this does not imply that the transversal is analytic (or locally algebraic). 
Assume that the continuation is given by $d_{1}d_{2}$, such that $S_{1} \rightarrow S_{2}$
through reflection with respect to a reflection axes in $\pi$. Transversals can be seen as 
$\bot \pi$. On the hyperboloid, we can always represent proper movements as reflection 
$\sim$ a normal model. Movement can be continued to euclidean metric, also the reflection axes, 
but orthogonals (corresponding to the transversal) are not necessarily analytic.

Note that in the problem of localisation, the condition of surjectivity is not necessary. 
It is sufficient with existence of $\gamma$, such that $f(\zeta) \rightarrow \gamma \rightarrow V$ continuous. 
Assume $\forall \gamma \in (J)$ we have existence of $R$ hypoelliptic with $R \gamma_{2}=\gamma$. 
Assume $P^{*}=R$ in $\mathcal{D}_{L^{2}}'$. 
When $E$ is very regular, we have that $E \gamma=\gamma_{2}$ hypoelliptic and $E \sim \delta_{0} - \eta$ 
implies $\gamma=\gamma_{2} + \eta * \gamma$ gives an approximative solution.
In this article, the mixed problem is dealt with as $F(\gamma_{1} \rightarrow \gamma_{2})(\zeta) 
=f(\zeta)$ and in the converse direction over $\gamma_{2} \rightarrow \gamma_{1}$, the mixed
problem is already extensively dealt with (for instance \cite{Ikawa85},\cite{Ikawa88}).

Assume $f_{1} \sim f + \delta$. This equivalence can be used in connection with Cousin's model 
of monotropy. For isolated singularities we have $f_{1}(\xi)=f(\xi) + \delta=f(\xi + \epsilon)$. 
We do no have in this case a normal covering in the sense of Weyl (\cite{Weyl09}), but 
using a reduced complement, if the singular points are $\notin \mbox{ supp }\mu$,
we can assume $f \bot \mu$ iff $f-\delta \bot \mu$. In connection with double transform, 
$P(D) \delta \rightarrow P(\xi) \rightarrow P(\xi^{*})$,
when $f(\zeta) + \delta(\zeta)=F(r_{T}' \gamma)(\zeta)$, where $F$ is linear in $\gamma$. 
Given $\gamma$ analytic, 
we can locally write $F \sim 1/Q$ for a polynomial $Q$, in this case $F(\gamma)$ is constant.
Given $\int g d \mu=0$ with $g$ algebraic, we have that the intersection is of measure zero. 
When $g$ is not algebraic, we are discussing one-sidedness. $F(g_{+})=G(g_{-})$
where the regularity is preserved.

\subsection{Continuation of the boundary}
Assume the boundary $\Gamma$ extended with mirrors (congruent to tangents) to 
$\tilde{\Gamma}$. Starting from a boundary point in $\Gamma$ selected as origo, 
the problem is to reach the others through a chain of broken rays. If also the chain 
is regarded as $\tilde{\Gamma}$ then we reach in this manner a subset of the boundary of 
positive measure (we assume the symbol $f=0$ over the chain). Assume that $d \mu$ a boundary measure 
(\cite{Garding64}) and denote with $\tilde{d \mu}$ the measure corresponding to the extended 
boundary. Thus $\tilde{\Gamma} \rightarrow \{ \mid z \mid=1 \}$

Consider now $A(x,y)=A(x(1,\frac{y}{x}))$. If $(1,\frac{y}{x}) \in \Omega$, where $\Omega$ 
is the domain for analyticity for $A$, we have that $A$ is analytic on $1 + \mid \frac{y}{x} \mid^{2} \leq 1$, 
In particular, if $f$ is analytic on $(1,\frac{y}{x})$ with $\mid \frac{y}{x} \mid < 1$ and 
on $(x,0)$ where $\mid x \mid < 1$, we have that $f$ is analytic on $\mid x \mid^{2} + \mid y \mid^{2} \leq 1$.
The symmetry condition $F(r_{T}' x,y)=F(x,r_{T}' y)$ means that the zero set is symmetric in 
$(x,y)$. If $F(x,y)=0$ we have that $F=0$ in $r_{T}' x$, for $y$ fixed, that is on a cone 
relative $r_{T}'$ (relative homogeneity).
The condition $y_T / x_{T}=y/x$ can be written $F(r_{T}' x,r_{T}' y)=F(x,y)$.
Assume $P$ hyperbolic and $Q$ partially hypoelliptic with $\Delta(P)=\Delta(Q)$, 
consider for a parmetrix $E_{1}$ to $P$, $\tau P v =P$, that is $vE_{1} = E_{1}$ over $\Delta$. 
If $E_{2}$ is parametrix to $Q$ and if we consider the parametrices as Fredholm operators, 
with the difference that $\mbox{ ker }E_{2}$ can be reduced to a trivial space through iteration.

Assume $\Delta^{\bot}$ is given by $S_{j}$, continued to $\tilde{S}_{j}$ through 
$M(e^{\psi} d \gamma) \rightarrow e^{\phi} M(d \gamma)$ continuous, 
where $\phi \in L^{1}$. If we have $\phi>0$, we have trivial first surfaces, 
corresponding a reduced operator. 
We are assuming $\tilde{S}_{j}$ can be continuously deformed to $S_{j}$ using $\phi$. 
If we assume the symmetry condition for $\phi$, we can
assume the deformation independent of parabolic movements.

If $F(\gamma)$ is of type (A), which is the case where $\gamma=\gamma_{1}$, we can chose 
$F(\gamma)=P(D) \tilde{f}_{0}$, where $\tilde{f}_{0}$ is very regular and the compatibility condition
according to Kiselman  (\cite{Kiselman65}) means that we have existence of $Q_{j} \bot P$. Assume $\phi_{j}$ corresponds to 
$\gamma_{2}^{j}$, such that in this case $\phi_{j} \rightarrow 0$, when $j \rightarrow \infty$. 
Note the difference between the zero space and constant surface, $\phi_{1}+\phi_{2}=const$. 
If we describe the first surfaces as $V_{1} \cup V_{2}$ where $V_{j}$ is connected, 
this corresponds to a multiply connected boundary, where every $\phi_{j}$
is assumed such that the transversal is locally represented as a polynomial 
(at least analytic). 

Assume $\mu$ a reduced measure of bounded variation and $\mu_{2}=P_{1} \ldots P_{N} \mu$, 
we have in this case not approximations on the form of orbits. In this case, even if we 
for monotropic functions do not have have presence of a normal covering in the sense of Weyl, 
we do not have problems with orbits, given that we start with a reduced measure in the 
representation of the boundary. 
Note the moment problem, if we assume $\alpha$ reduced, we then have $\alpha(g)=0$ 
implies $g=0$ on the domain $\Omega$, that is given $\alpha \in \mathcal{D}_{L^{1}}'$ is 
$\Omega$ a strict carrier to $\alpha$. Since $g \equiv 0$ on $\Omega$ 
or $\alpha(\Omega)=0$, if $g$ is a polynomial and if $\alpha$ reduced, we must $\Omega$ is trivial. 
Finally, if the condition is $\phi \in L^{1}$ and we consider $\psi \in H(V)$ such that 
$\psi \sim_{m} \phi$, we do not longer have a normal covering,
but according to the moment problem, we can uniformly approximate $C \cap L^{1}$ with 
$H(V)$, over a strict carrier to the measure $\alpha$. If we only assume the measure 
$\alpha$ of bounded variation, but instead we chose $g$ as polynomial, we can assume '
$\Omega$ has $\alpha-$ measure zero.

Given a measure of bounded variation and positive definite $d \mu$, such that 
$\tilde{d \mu}=d \mu + d \mu_{0}$, where $d \mu_{0}$ is assumed with point support and 
$\tilde{d \mu}$ holomorphic, we have a continuation according to Cousin (\cite{Dahn13}, compare also \cite{Garding64}).
When $\tilde{ d \mu}$ locally reduced, we have $\int_{\gamma_{T} - \gamma_{0}} \tilde{d \mu}=0$ 
implies $\gamma_{T}=\gamma_{0}$.
When $\tilde{d \mu}$ is not reduced, there is the centre case among the possible approximations, 
that is orbits around a singular point. 
Given $g$ a regular approximation and 
$\tilde{ d \mu}$ with support on first surface (the boundary), assume $\int g \tilde{d \mu}=0$ 
under the approximations. Given $\tilde{d \mu}$ reduced, we can assume $g=0$ outside a compact set. 
It is clear that $L$ (transversals) are included in the support for $g$, in this manner we can regard $L$ as a strict carrier for limits. 
The proposition is that given the transversal as a strict carrier for limits, the transversal can 
be represented through polynomial locally.  

The boundary is assumed very regular in the sense Parreau (\cite{Parreau51}). We assume all 
singularities on first surfaces and that all first surfaces can be reached. 
In the context of the moment problem, we can allow monotropical functions, 
that is $g$ can be continued to continuous functions.

Assume $M=\tilde{d \mu}/d z$ and $\frac{d G}{d z}=g$, we then have if $g$ regular that $dG$ 
is a closed form. 
Assume $\Sigma=\{ f(\zeta)=const \}$ and $J U = V J$ according to the invariance principle
and $S=\{ F(\gamma_{1})=const \}$ and $\tilde{S}=\{ F(\gamma_{1} \rightarrow \gamma_{2})=const \}$.
If $\varphi$ is a continuous function in a neighbourhood of $(x^{0},y^{0})$ and $\sigma$ 
a characteristic surface through this point, we can then write $\sigma \subset \{ \varphi(x,y) >0 \} \backslash (0,0)$ 
(pseudo convexity) and $\sigma=\{ x_{n}=P(x_{1},\ldots,x_{n-1}) \}$ (strict pseudo convexity (\cite{Oka60})). 
We define $U \gamma_{1}=V \gamma_{2}$ in a neighbourhood of $\Delta=\{ U=V=I \}$ and assume 
$F(U \gamma_{1}) \rightarrow const$. For all singularities that can be reached in this manner, 
given $U$ is a fixed movement, we can chose $\gamma_{1} \rightarrow \gamma_{2}$ as analytic 
(polynomial) .
Note that, if $\Delta$ is a plane $x_{j}=0$, we have that $e^{v}$ is constant in some direction, 
why the function is not in $L^{1}$ on an unbounded domain.

Given $0 = \int_{V} g d \mu$ we have that $g \equiv 0$ on $V$ or $\mu(V)=0$ (measure zero).
Given $g$ polynomial (exponential of a polynomial) we have from a result by Hurwitz (\cite{Dahn13})
that the measure for all singularities is zero. In the transposed case $0 = \int M d g$ 
implies $M=0$ or $g(V)=0$. If for instance $M$ is locally 1-1, we have that $x=x_{0}$.

\section{The moment problem}
Note the moment problem, assume $E_{0}$ a set in $C$, the problem can be solved if we existence of $\alpha(x)$ of bounded variation, such that $\int g(x) d \alpha(x)=0$
for $g \in E_{0}$, we have $\int g(x) d \alpha(x)=0$ for $g \in C$ (\cite{Riesz56}). Note that 
when the orthogonal (transversal) is algebraic, this corresponds to a measure zero complement .

Concerning mixed problems, assume $(Au,v)=(u,A^{*}v)$, such that $A^{*}v=0$, that is 
$v \bot Au=f$. Given that $Bu=v$, we have thus  $A \bot B$ or $A^{*}B \bot I$. 
Note $R(B)=D(A)^{\bot}$ and $R(A)=D(B)^{\bot}$. If 
$R(A)^{\bot}=R(B)$ we have that $B$ has closed range, even when $A$ does not.
If $B$ is algebraic in $\mathbf{R}^{n}$, we have that given $v$, we have existence of $u$ such that
$Bu=v$. The domain for $\bot$ is then $u$, such that $u \in D(A) \cap D(B)$ that is closed
when $D(A)$ is closed. Note if $A$ is a polynomial, there are complex polynomials $B^{*}$, 
such that the symbol to $B^{*}A \equiv 0$ (\cite{Kiselman65}). If we consider $B^{*}$ as annihilator, $(B^{*}A u,u)=0$ 
for all $u \in D(A)$. When $A$ is considered in $X'$, we have $B^{*}A = 0$, that is $B^{*}$ is an 
annihilator for $A$. 
Compare with the moment problem, where a sufficient condition, given finite order singularities, for existence of $B^{*}$ is an analytic 
representation of $A$ (Parseval). 
Consider now $A$ hyperbolic and $B^{N}$ hypoelliptic, it is sufficient to consider $u$ with $A^{N}u=f_{N}$ and $B^{N} v=0$ for $v \bot f$.
Note that if $E_{N}$ is a parametrix to $B^{N}$, it has (modulo $C^{\infty}$) a trivial kernel, that is
if $v'=E_{N}v$, $f_{N} \bot v$, we have $v'=0$.

\newtheorem*{res14}{The mixed model}
\begin{res14}
 Assume $\gamma \rightarrow \gamma_{\delta}$ is a continuation of $\gamma$ and $F_{\delta}$ is constructed,
 so that $F_{\delta}(\gamma_{\delta})=f$ and respects $f \bot v$, where $v$ is continuous. Assume ${}^t F_{\delta}$ has a trivial
 kernel, then the continuation gives a global mixed model.
\end{res14}

Note that $\delta_{0} - C^{\infty}$ maps $\mathcal{D}' \rightarrow \mathcal{D}^{' F}$,
that is if the boundary is defined modulo regularizing action instead of modulo $H$, 
we do not longer have that $r_{T}$ is injective in $H$, however it is injective in $L^{1}$. 
We will use the moment problem (and when possible the transmission problem (\cite{Schechter60})) to solve the mixed problem. In the two mirror model
we assume that a hyperbolic operator $A$, is reflected through the boundary into an operator $A^{*}$ 
(geometrical optics) and in the same manner for $B$ partially hypoelliptic. Existence of a continuous 
mapping between the respective boundary points, implies an abstract transmission property. 
Note that only one of the systems has to be invertible.

Assume the condition $B \gamma_{1}=0$ on $\Sigma_{i}$ corresponds to $\Delta(f)$ (lineality), then there is a $B'$ such that 
$B \gamma_{1}=0$ iff $B' \gamma_{2}=0$ on $\Sigma_{i}$. Thus, if $\gamma_{2}$ is seen as a 
continuation of $\gamma_{1}$ according to $F(\gamma_{1} \rightarrow \gamma_{2})$
with $\gamma_{1} \bot \gamma_{2}$, we have a ``global model''. We have $< F(\gamma_{1}), v>=0$ and
${}^{t} F(v) \bot \gamma_{1}$ and ${}^{t} F(v)=0$ implies $v \bot F(\gamma_{1})$,
where $F$ is a parametrix (localization) to $\gamma_{1} \rightarrow \gamma_{2}$. Further, $< F(\gamma_{2}),v>=0$ 
and $v \in \mbox{ ker } {}^{t} F=\{ 0 \}$, which implies $v=0$. Assume $J_{\delta} : \gamma_{1} \rightarrow \gamma_{2}$, 
when $\delta \rightarrow 0$ and $J_{\delta} \gamma_{1} \bot \gamma_{1}$, further that $({}^{t} J_{\delta} F) \rightarrow E$ 
with $\mbox{ ker }E=\{ 0 \}$. Define $F_{\delta} \sim ({}^{t} J_{\delta} F)$ with 
$< \gamma_{1}, {}^{t} F_{\delta} v > \rightarrow 0$, as $\delta \rightarrow 0$. According to the above, 
this implies $v=0$, which implies $F_{\delta}=id$. Thus, if $F_{\delta}(\gamma_{1})=f$ with $F_{\delta}=id$, 
then this implies that $f$ hyperbolic and when $F_{\delta}(\gamma_{2})=f$ with $F_{\delta}=id$, then this 
implies that $f$ hypoelliptic.

Generally, for a radical geometric ideal, it is sufficient to give the boundary condition on derivatives to the symbol. 
Consider otherwise the problem for $M(f)$ (arithmetic mean).  
The boundary condition, when it is radical, can be given for the derivative, but does not necessarily define
the domain. The radical boundary condition, does not imply a radical ideal. 
Consider $J_{\delta} : \gamma_{1} \rightarrow \gamma_{2}$, such that $\lim_{\delta \rightarrow 0} J_{\delta} \gamma_{1} \rightarrow $
a hypoelliptic symbol. Write $J_{\delta} \gamma_{1}=\gamma_{\delta}$, we then have $F(J_{\delta} \gamma_{1})=
\big[ F , E_{\delta} \big](\gamma_{2})$, where $E_{\delta} \sim e^{v_{\delta}}$, where $v_{\delta}=0$
on invariant sets. On the other hand, $F \gamma_{2} = 1$, as $v_{\delta}=0$ implies $\mbox{ ker }F = \{ 0 \}$.
Assume $F(\gamma) \rightarrow \gamma \rightarrow \zeta$. When $U \gamma$ is analytic, we have $\tilde{U} \zeta$ 
continuous. Note that $\{ \zeta \quad U \gamma_{1} = \gamma_{1} \}= \{ \zeta \quad U \gamma_{1}^{2}= \gamma_{1}^{2} \}$
and $\{ \zeta \quad V \gamma_{2}^{2} = \gamma_{2}^{2} \} \subset \{ \zeta \quad V \gamma_{2} = \gamma_{2} \}$. As $J_{\delta} U \gamma_{1}^{2} = V J_{\delta} \gamma_{1}^{2} = V \eta_{2}$
where $\eta_{2}=J_{\delta} \gamma_{1}^{2}$ and $\eta_{2} \in (I_{Phe})$. 

Consider the continuation $S_{j} \rightarrow \tilde{S}_{j}$ simply connected. 
Given Oka's condition for $\tilde{S}_{j}$ and the continuation of the symbol, the transversal can be defined
as analytic (algebraic) even for $\tilde{S}_{j}$. Concerning the set $\{ f < \lambda \}=V$,
if we have $S_{j} \subset \subset V$ and $S_{j}$ is a bounded set, we do not necessarily have the same
for $\tilde{S}_{j}$. The proposition for the continuation implies a downward bounded symbol.

According to Weyl (\cite{Weyl09}), we have for a normal covering, that closed curves corresponds to 
closed curves. According to (\cite{Poincare87}) (Ch. XXVII), we have modulo monotropy, that closed curves 
correspond to open curves with points in common with its closed correspondent. 
We can use the moment problem, given $\int_{(\Gamma)} g d \alpha=0$ for a reduced 
measure of bounded variation,
for $0 \neq g \in E_{0}$ implies the same relation for $g \in C$. Thus, in the plane we can use a
reduced measure (modulo removable sets) to solve the problem.

Assume $(I)=(\mbox{ ker }h)$ and $d h(f)=g(z)d z$, where $g \sim_{m} g_{1}$ analytic. Assume 
$g-g_{1}$ algebraic and $g_{1}$ regular over closed contours. We can assume one-sided regularity for 
$g$. For instance, $f(1/x)g(x)$
where $f$ is bounded close to the $\infty$ and $g$ is bounded close to $0$.
Assume a two-sided limit and $F_{+}(g_{+}) - F_{-}(g_{-}) \rightarrow 0$ at the boundary. Assume 
$g_{-} \sim g_{+}^{*}$ and $F_{-}(g_{-}) \sim F_{+}^{*}(g_{+})$, we then have $F_{+} - F_{+}^{*}(g_{+}) \rightarrow 0$ 
and through Radon-Nikodym's theorem, we have $F_{+}(g_{+})=F_{-}(e^{v} g_{-})$, where $v \in L^{1}$.

The Schr\"odinger operators (\cite{Riesz46}) give a global model, that is not normal. 
If we select $d \mu$ very regular, that is hypoelliptic in $L^{2}$,
the representation of the symbol is locally 1-1, why we have two-sided limits. 
Compare with $\int g (d \mu - d \mu_{0})=0$, where $d \mu_{0}$ has 
point support. Then $d \mu = v d z$ and $< g,d \mu>=<g,v>$ and $v=v_{1} + v_{0}$, we can assume 
$v_{0} \sim \delta / g$. Note that $V^{\bot} \gamma \rightarrow x_{0}$ does not imply $V \gamma \rightarrow x_{0}$. 
Consider $V^{\bot} \rightarrow x^{*}$ continuous and $V \rightarrow x$, further that $x_{0}$ is a
joint point. For instance, if $\widehat{f}(x^{*}) \rightarrow 1$, we have simultaneously 
$f(x) \rightarrow \delta$. We can then have $x^{*} \rightarrow x_{0}$ without  simultaneously $x \rightarrow x_{0}$

Concerning pseudo vectors, consider a neighbourhood of a point on $H_{m}$. An infinitesimal displacement 
can be performed as parabolic, elliptic or hyperbolic. The normal is considered relative an axes 
for invariance, that is we consider the normal as independent of the neighbourhood. Assume for $x,y$ real,
$\eta(x)=y(x)/x = e^{\phi(x)}$. We have then three possibilities, $\phi=0$ parabolic, $\phi<0$ hyperbolic and $\phi>0$ elliptic.
The three possibilities induce three possible orientations for the normal. In particular in the parabolic case,
if the normal is dependent on scaling parameter $h$, as $h \rightarrow 0$ or $h \rightarrow \infty$, we have
a `` scaling orientation'' for the normal. In the case where the transversal is a strict carrier or a carrier, 
the limit is not dependent of choice of neighbourhood.

Form the normal in a point $\zeta_{0}$ at the boundary $\Gamma=\{ \zeta \quad F(\gamma)(\zeta) = const. \}$
and consider $\gamma(\zeta_{j} + h) - \gamma(\zeta_{j} - h)$, where $h$ is a scalar. Obviously, if
$\gamma$ is symmetric with respect to $\zeta_{0}$ over the segment, we have that $\gamma$ is absolute
continuous, as $h \rightarrow 0$. In this case the definition of the normal, does not depend on $h$.
Consider $G(h)=\int N(\widetilde{\zeta}) \gamma'(\widetilde{\zeta}) d \widetilde{\zeta}$,
where for instance $\widetilde{\zeta}=(\zeta_{1},\ldots,\zeta_{j}+h,\ldots,\zeta_{n})$. When $N$ is polynomial
over the segment $I_{h}$, we have when $\gamma$ is absolute continuous and $G(h) \equiv 0$, that $\gamma=const.$
over $I_{h}$. Further, $N \bot d \gamma$, independent of $h \rightarrow 0$. Consider now $N \in \mathcal{D}_{L^{1}}'$,
summable with respect to $d F$, of bounded variation, we then have that $G$ is absolute continuous with respect to $dF$,
that is $F(I_{h}) \rightarrow 0$, as $h \rightarrow 0$, implies $G(h) \rightarrow 0$. When $G \equiv 0$,
we have that $N \bot d F$. Thus, when $N \in L^{1}(dF)$, we have if $G$ absolute continuous with respect to $dF$,
that $N \bot dF$, as $h \rightarrow 0$. When $N$ is dependent on $h$, we may have $G \nrightarrow 0$, as $h \rightarrow 0$.
Further, when $N$ has a strict carrier, then $N$ is independent of $h$ for large $h$. This ``orientation''
of $N$ is represented using regularity conditions.

\bibliographystyle{amsplain}
\bibliography{mixed}

\providecommand{\bysame}{\leavevmode\hbox to3em{\hrulefill}\thinspace}
\providecommand{\MR}{\relax\ifhmode\unskip\space\fi MR }
\providecommand{\MRhref}[2]{%
  \href{http://www.ams.org/mathscinet-getitem?mr=#1}{#2}
}
\providecommand{\href}[2]{#2}
\begin{thebibliography}{10}

\bibitem{AhlforsSario60}
L.~Sario~L. Ahlfors, \emph{{Riemann Surfaces}}, Princeton University Press,
  1960.

\bibitem{Bendixsson01}
I.~Bendixsson, \emph{{Sur Les Courbes D\'{e}finies Par Les \'{E}quations
  Diff\'{e}rentielles.}}, Acta Math. \textbf{24} (1901).

\bibitem{Dahn13}
T.~Dahn, \emph{{Some remarks on Treve's conjecture}}, ArXiv \textbf{1303.3802}
  (2013).

\bibitem{Riesz56}
B.~Sz.-Nagy {F. Riesz}, \emph{{Functional Analysis.}}, Dover Publications Inc.,
  1956.

\bibitem{Garding87}
L.~G{\aa}rding, \emph{{Singularities in Linear Wave Propagation.}}, Springer
  Lecture Notes in Mathematics, 1241, 1987.

\bibitem{Garding64}
L.~G{\aa}rding and L.~H\"{o}rmander, \emph{{Strongly Subharmonic Functions.}},
  Math. Scand. \textbf{15} (1964).

\bibitem{Riesz15}
G.H. Hardy and M.~Riesz, \emph{{The General Theory of Dirichlet Series.}},
  vol.~18, Cambridge Tracts in Mathematics and Mathematical Physics, 1915.

\bibitem{Ikawa85}
M.~Ikawa, \emph{{Sur la d\'{e}croissance d'\'{e}nergie locale du probl\`{e}me
  ext\'{e}rieur avec plusieurs $(n \geq 3)$ obstacles strictement convexes}},
  S\'{e}minaire de th\'{e}orie spectrale et g\'{e}om\'{e}trie \textbf{4}
  (1985-1986).

\bibitem{Ikawa88}
\bysame, \emph{{Decay of solutions of the wave equation in the exterior of
  several convex bodies}}, Annales de l'institut Fourier \textbf{38} (1988),
  no.~2.

\bibitem{Janc89}
B.~Jancewicz, \emph{{Multivectors and Clifford Algebra in Electrodynamics.}},
  World Scientific Publishing Co Pte Ltd, 1989.

\bibitem{Julia23}
G.~Julia, \emph{{Lecons sur les fonctions uniformes \`{a} point singulier
  essentiel isol\'{e}}}, Gauthier-Villars, Paris, 1923.

\bibitem{Julia24}
\bysame, \emph{{Sur quelques applications de la repr\'{e}sentation conforme
  \`{a} la r\'{e}solution d'\'{e}quations fonctionelles}}, Bulletin de al
  Soci\'{e}t\'{e} Math\'{e}matique de France \textbf{52} (1924).

\bibitem{Kiselman65}
C.O. Kiselman, \emph{{Existence and approximation theorems for solutions of
  complex analogues of boundary problems}}, Arkiv f\"{o}r Matematik \textbf{6}
  (1965), no.~11.

\bibitem{Schechter60}
Schechter M., \emph{{A generalization of the problem of transmission}}, Annali
  della Scuola Normale Superiore di Pisa- Classe di Scienze \textbf{14} (1960),
  no.~3.

\bibitem{Martineau}
A.~Martineau, \emph{{Sur Les Fonctionelles Analytiques et la Transformation de
  Fourier-Borel}}, Journal d'Analyse Math\'{e}matique, Jerusalem \textbf{XI}
  (1963).

\bibitem{Nishino68}
T.~Nishino, \emph{{Nouv\`{e}lles Recherches sur Les Fonctions Entires de
  Plusieurs Variables Complexes.}}, Journal Math. Kyoto Univ. (1968).

\bibitem{Nishino75}
\bysame, \emph{{Nouvelles recherches sur les fonctions enti\`{e}res de
  plusieurs variables complexes (V) Fonctions qui se r\'{e}duisent aux
  polynomes}}, J. Math. Kyoto Univ. \textbf{15} (1975), no.~3.

\bibitem{Oka60}
K.~Oka, \emph{{Sur Les Fonctiones Analytiques de Plusieurs Variables}},
  Hermann, 1960.

\bibitem{Parreau51}
M.~Parreau, \emph{{Sur les moyennes des fonctions harmoniques et analytiques et
  la classification des surfaces de Riemann}}, Annales de l'institut Fourier
  \textbf{3} (1951).

\bibitem{Poincare87}
H.~Poincar\'{e}, \emph{{Les M\'{e}thodes Nouvelles de la M\'{e}canique
  C\'{e}leste.}}, Gauthier-Villars, 1987.

\bibitem{Riesz43}
M.~Riesz, \emph{{En {\aa}sk{\aa}dlig bild av den icke-euklidiska geometrien.
  Geometriska str\"{o}vt{\aa}g inom relativitetsteorien.}}, {Lunds Universitets
  {\AA}rsskrift N.F. Avd. 2, Bd. 38, Nr. 9 }, vol. Lund, C.W.K. Gleerup, 1943.

\bibitem{Riesz46}
\bysame, \emph{{Sur certaines notions fondamentals en th\'{e}orie quantiques
  relativiste.}}, 10e Congr. des Math. Scand. Copenhagen (1946).

\bibitem{Schwartz66}
L.~Schwartz, \emph{{Th\'{e}orie Des Distributions.}}, Hermann, 1966.

\bibitem{Vuorinen07}
M.~Vuorinen, G.D. Anderson, and M.K. Vamanamurthy, \emph{{Generalized convexity
  and inequalities}}, arXiv:math \textbf{0701262v1} (2007).

\bibitem{Weyl09}
H.~Weyl, \emph{{The concept of a Riemann Surface.}}, Dover Publications Inc.,
  2009.

\end{thebibliography}

\end{document}